\documentclass[aos,MSNbibl,seceqn,citesort,noautosecdot,rotating,dvips]{arximspdf}
\usepackage{mathbh}
\usepackage{dcolumn}
\usepackage{graphicx}

\doi{10.1214/11-AOS918}
\volume{39}
\issue{5}
\pubyear{2011}
\firstpage{2766}
\lastpage{2794}

\makeatletter

\newcolumntype{d}[1]{D{.}{.}{#1}}

\newcommand{\Tr}{\operatorname{Tr}}
\newcommand{\Diag}{\operatorname{Diag}}
\newcommand{\ess}{\operatorname{ess}}
\newcommand{\sign}{\operatorname{sign}}
\newcommand{\Span}{\operatorname{span}}

\newcommand{\YY}{\bolds{\mathsf{Y}}}

\newtheorem{theorem}{Theorem}[section]

\newcommand{\C}[1]{\mathcal{#1}}
\newcommand{\wt}[1]{\tilde{#1}}
\newcommand{\wh}[1]{\hat{#1}}

\newcommand{\undc}[2]{{#1}_{#2}}

\newcommand{\eqdef}{\triangleq}

\newcommand{\cC}{\Theta}
\newcommand{\cD}{\mathcal{D}}
\newcommand{\E}{\mathbb{E}}
\newcommand{\cF}{\mathcal{F}}

\newcommand{\M}{\mathcal{M}}
\newcommand{\N}{\mathcal{N}}
\renewcommand{\P}{\mathbb{P}}
\newcommand{\R}{\mathbb{R}}
\newcommand{\cS}{\mathcal{S}}
\newcommand{\bcR}{B} 
\newcommand{\cR}{\tilde{B}} 
\newcommand{\X}{\mathcal{X}}
\newcommand{\Z}{\mathcal{Z}}

\newcommand{\hth}{\hat{\theta}}

\newcommand{\kap}{\kappa}
\newcommand{\lam}{\lambda}

\newcommand{\logeps}{\log(\varepsilon^{-1})}

\newcommand{\ra}{\rightarrow}

\newcommand{\argmin}{\arg\min}

\newcommand{\vp}{\varphi}
\renewcommand{\th}{\theta}
\newcommand{\hf}{\hat{f}}
\newcommand{\chr}{\check{r}}

\newcommand{\hfols}{\hf^{(\mathrm{ols})}}
\newcommand{\therm}{\hat{\th}^{(\mathrm{erm})}}
\newcommand{\hferm}{\hf^{(\mathrm{erm})}}

\newcommand{\thrid}{\tilde{\th}} 
\newcommand{\frid}{\tilde{f}} 
\newcommand{\freg}{f^{(\mathrm{reg})}}
\newcommand{\thrlam}{\hth^{(\mathrm{ridge})}} 
\newcommand{\hfrlam}{\hf^{(\mathrm{ridge})}} 
\newcommand{\Flin}{\mathcal{F}_{\mathrm{lin}}}
\newcommand{\XX}{\bolds{\mathsf{X}}}

\newcommand{\V}[1]{\overline{#1}}
\newcommand{\bL}{\V{L}}

\makeatother

\begin{document}
\begin{frontmatter}

\title{Robust linear least squares regression}
\runtitle{Robust linear least squares regression}

\begin{aug}
\author[A]{\fnms{Jean-Yves} \snm{Audibert}\ead[label=e1]{audibert@imagine.enpc.fr}}
\and
\author[B]{\fnms{Olivier} \snm{Catoni}\corref{}\ead[label=e2]{olivier.catoni@ens.fr}}
\runauthor{J.-Y. Audibert and O. Catoni}
\affiliation{Universit\'e Paris-Est and CNRS/\'Ecole Normale Sup\'
erieure/INRIA\break and
CNRS/\'Ecole Normale Sup\'erieure and INRIA}
\address[A]{Universit\'e Paris-Est\\
LIGM, Imagine\\
6 Avenue Blaise Pascal\\
77455 Marne-la-Vall\'ee\\
France\\
and\\
CNRS/\'Ecole Normale Sup\'erieure/INRIA\\
LIENS, Sierra---UMR 8548\\
23 Avenue d'Italie\\
75214 Paris Cedex 13\\
France\\
\printead{e1}} 
\address[B]{D\'epartement de Math\'ematiques\\
\quad et Applications\\
\'Ecole Normale Sup\'erieure\\
CNRS---UMR 8553\\
45 Rue d'Ulm\\
75230 Paris Cedex 05\\
France\\
and\\
INRIA Paris-Rocquencourt---CLASSIC team\\
\printead{e2}}
\end{aug}

\received{\smonth{2} \syear{2009}}
\revised{\smonth{8} \syear{2011}}

%
\begin{abstract}
We consider the problem of robustly predicting as well as the best
linear combination of $d$ given functions in least squares regression,
and variants of this problem including constraints on the parameters of
the linear combination. For the ridge estimator and the ordinary least
squares estimator, and their variants, we provide new risk bounds of
order $d/n$ without logarithmic factor unlike some standard results,
where $n$ is the size of the training data. We also provide a~new
estimator with better deviations in the presence of heavy-tailed noise.
It is based on truncating differences of losses in a min--max framework
and satisfies a $d/n$ risk bound both in expectation and in deviations.
The key common surprising factor of these results is the absence of
exponential moment condition on the output distribution while achieving
exponential deviations. All risk bounds are obtained through a
PAC-Bayesian analysis on truncated differences of losses. Experimental
results strongly back up our truncated min--max estimator.
\end{abstract}

%

\end{frontmatter}

\section{\texorpdfstring{Introduction.}{Introduction}}

\subsection*{\texorpdfstring{Our statistical task.}{Our statistical
task}}

Let $Z_1=(X_1,Y_1),\ldots,Z_n=(X_n,Y_n)$ be $n\ge2$ pairs of input--output
and assume that each pair has been independently drawn from the same
unknown distribution $P$. Let
$\X$ denote the input space and let the output space be the set of
real numbers $\R$, so that $P$ is a probability distribution on the
product space
$\Z\eqdef\X\times\R$.
The target of learning algorithms is to predict the output $Y$
associated with an input $X$
for pairs $Z=(X,Y)$ drawn from the distribution $P$.
The quality of a (prediction) function
$f\dvtx\X\rightarrow\R$ is measured by the least squares \textit{risk}:
\[
R(f) \eqdef\undc{\E}{Z\sim P} \{ [Y-f(X)]^2
\} .
\]
Through the paper, we assume that the output and all the prediction
functions we consider are square integrable.
Let $\cC$ be a closed convex set of~$\R^d$, and $\vp_1,\ldots,\vp
_d$ be $d$ prediction functions.
Consider the regression model
\[
\cF= \Biggl\{ f_\th=\sum_{j=1}^d \th_j \vp_j ; (\th_1,\ldots,\th
_d) \in\cC\Biggr\}.
\]
The best function $f^*$ in $\cF$ is defined by
\[
f^* =\sum_{j=1}^d \th^*_j \vp_j\in\mathop{\argmin}_{f\in\cF} R(f).
\]
Such a function always exists but is not necessarily unique. Besides,
it is unknown since the probability generating the data is unknown.

We will study the problem of predicting (at least) as well as function~%
$f^*$. In other words, we want to deduce from the observations
$Z_1,\ldots,Z_n$ a~function~$\hf$ having with high probability a risk
bounded by the minimal risk $R(f^*)$ on $\cF$ plus a small remainder
term, which is typically of order $d/n$ up to a~possible logarithmic factor.
Except in particular settings (e.g., $\cC$ is a~simplex and $d\ge
\sqrt{n}$),
it is known that the convergence rate $d/n$ cannot be improved in a
minimax sense (see \cite{Tsy03} and \cite{Yan01b} for related results).

More formally, the target of the paper is to develop estimators $\hf$
for which the excess risk is controlled \textit{in deviations}, that is,
such that
for an appropriate constant $\kap>0$, for any $\varepsilon>0$, with
probability at least $1-\varepsilon$,
%
%
\begin{equation} \label{eqdevtarget}
R(\hf) - R(f^*) \le\frac{\kap[d+\logeps]}{n}.
\end{equation}
Note that by integrating the deviations [using the identity
$\E( W ) =\break \int_0^{+\infty} \P(W > t) \,dt$
which holds true for any non-negative random variable~$W$],
inequality~(\ref{eqdevtarget}) implies
%
%
\begin{equation} \label{eqexptarget}
\E R(\hf) - R(f^*) \le\frac{\kap(d+1)}{n}.
\end{equation}
In this work, we do not assume that the function
\[
\freg\dvtx x \mapsto\E[ Y | X=x],
\]
which minimizes the risk $R$ among all possible measurable functions,
belongs to the model $\cF$. So we might
have $f^* \neq\freg$ and in this case, bounds of the form
%
%
\begin{equation} \label{eqnottarget}
\E R(\hf) - R\bigl(\freg\bigr) \le C \bigl[R(f^*)-R\bigl(\freg\bigr) \bigr] + \kap
\frac{d}{n}
\end{equation}
with a constant $C$ larger than $1$, do not even ensure that $\E R(\hf
)$ tends to~$R(f^*)$ when $n$ goes to infinity.
These kinds of bounds with $C>1$ have been developed to analyze
nonparametric estimators using linear approximation spaces,
in which case the dimension $d$ is a function of $n$ chosen so\vadjust{\goodbreak} that the
bias term $R(f^*)-R(\freg)$ has the order $d/n$ of the estimation term
(see \cite{Bar00,Gyo04,Sau10} and references within). Here we intend to assess the
generalization ability of the estimator even when the model is misspecified
[namely, when $R(f^*) > R(\freg)$].
Moreover, we do not assume either that $Y-\freg(X)$ and $X$ are
independent or
that $Y$ has a subexponential tail distribution: for the moment, we
just assume that $Y-f^*(X)$ admits a finite second-order moment in
order that the risk of $f^*$
is finite.

Several risk bounds with $C=1$ can be found in the literature.
A survey on these bounds is given in \cite{AudCat210}, Section 1.
Let us mention here the closest bound to what we are looking for.
From the work of Birg\'e and Massart \cite{BirMas98}, we may derive
the following risk bound for the empirical risk minimizer on a
$L^\infty$ ball
(see Appendix~B of \cite{AudCat210}).
\begin{theorem} \label{thbmnew}
Assume that $\cF$ has a diameter $H$ for $L^\infty$-norm, that is,
for any
$f_1,f_2$ in $\cF$, ${\sup_{x\in\X}} |f_1(x)-f_2(x)| \le H$ and there
exists a function
$f_0 \in\cF$ satisfying the exponential moment condition
%
%
\begin{equation} \label{eqexpmom}
\mbox{for any } x\in\X\qquad
\E\{ \exp[ A^{-1} \vert Y-f_0(X) \vert
] | X=x \} \le M
\end{equation}
for some positive constants $A$ and $M$. Let
\[
\cR= \inf_{\phi_1,\ldots,\phi_d} \sup_{\th\in\R^d - \{0\}}
\frac{\|{\sum_{j=1}^d \th_j \phi_j}\|_\infty^2}{\|\th\|_\infty^2},
\]
where the infimum is taken with respect to all possible orthonormal
bases of $\cF$ for
the dot product $(f_1,f_2) \mapsto\E[ f_1(X)f_2(X) ]$
(when the set $\cF$ admits no basis with exactly $d$ functions, we set
$\cR=+\infty$).
Then the empirical risk minimizer satisfies for any $\varepsilon>0$,
with probability at least $1-\varepsilon$,
\[
R\bigl( \hferm\bigr) - R(f^*) \le\kap(A^2 + H^2) \frac{d \log[2+(\cR
/n)\wedge(n/d)] +\logeps}{n},
\]
where $\kap$ is a positive constant depending only on $M$.
\end{theorem}

The theorem gives exponential deviation inequalities of order at worse
$d \log(n/ d)/{n}$ and, asymptotically, when $n$ goes to
infinity, of order $d/n$.
This work will provide similar results under weaker assumptions on the
output distribution.

\textit{Notation.}
When $\cC=\R^d$,
the function $f^*$ and the space $\cF$ will be written~%
$f^*_{\mathrm{lin}}$ and $\Flin$ to emphasize that
$\cF$ is the whole linear space spanned by $\vp_1,\ldots,\vp_d$:
\[
\Flin= \Span\{\vp_1,\ldots,\vp_d\} \quad\mbox{and}\quad
f^*_{\mathrm{lin}}\in
\mathop{\argmin}_{f\in\Flin} R(f).
\]
The Euclidean norm will simply be written as \mbox{$\|\cdot\|$}, and $\langle
\cdot,\cdot\rangle$
will be its associated inner product.
We will consider the vector valued function
\mbox{$\vp\dvtx\C{X} \rightarrow\mathbb{R}^d$} defined
by $\vp(X) = [ \vp_k(X) ]_{k=1}^d$, so that for any $\th
\in\Theta$, we have
\[
f_\th(X) = \langle\th, \vp(X) \rangle.\vadjust{\goodbreak}
\]
The Gram matrix is the $d\times d$-matrix $Q=\mathbb{E} [ \vp
(X) \vp(X)^T]$.
The empirical risk of a function $f$ is
$
r(f) = \frac{1}{n} \sum_{i=1}^n [ f(X_i) -
Y_i ]^2
$
and for $\lam\ge0$, the ridge regression estimator on $\cF$ is
defined by
$\hfrlam=f_{\thrlam}$ with
\[
\thrlam\in\mathop{\arg\min}_{\th\in\Theta} \{
r(f_{\theta}) + \lambda\Vert\theta\Vert^2 \},
\]
where $\lam$ is some non-negative real
parameter. In the case when $\lambda= 0$,
the ridge regression $\hfrlam$ is nothing
but the empirical risk minimizer $\hferm$.
Besides, the empirical risk minimizer when $\Theta=\R^d$
is also called the ordinary least squares estimator, and will be
denoted by
$\hfols$.

In the same way,\vspace*{1pt} we introduce the optimal ridge function
optimizing the expected ridge risk: $\frid=f_{\thrid}$ with
%
%
\begin{equation} \label{eqfrid}
\thrid\in\mathop{\arg\min}_{\theta\in\Theta} \{ R(f_{\theta})
+ \lambda\Vert\theta\Vert^2 \}.
\end{equation}
Finally, let $Q_{\lambda} = Q + \lambda I$ be the ridge regularization
of $Q$, where $I$ is the identity matrix.


\subsection*{\texorpdfstring{Why should we be interested in this task?}{Why should we be interested in this
task}}

There are four main reasons.
First, we intend to provide a nonasymptotic analysis of the parametric
linear least squares method.
Second, the task is central in nonparametric estimation for linear
approximation spaces (piecewise polynomials based
on a~regular partition, wavelet expansions, trigonometric
polynomials$\ldots$).

Third, it naturally arises in two-stage model selection. Precisely,
when facing the data, the statistician often has to choose
several models which are likely to be relevant for the task.
These models can be of similar structure (like embedded balls of
functional spaces) or,
on the contrary, of a very different nature (e.g., based on kernels,
splines, wavelets or on a parametric approach).
For each of these models, we assume that we have a learning scheme
which produces
a ``good'' prediction function in the sense that it predicts as well as
the best
function of the model up to some small additive term. Then the question is
to decide on how we use or combine/aggregate these schemes.
One possible answer is to split the data into two groups, use
the first group to train the prediction function associated with each
model, and
finally use the second group to build a prediction function which is as
good as (i) the
best of the previously learned prediction functions, (ii) the best
convex combination
of these functions or (iii) the best linear combination of these functions.
This point of view has been introduced by Nemirovski in \cite{Nem98}
and optimal rates of aggregation are given in \cite{Tsy03} and the
references within.
This paper focuses more on the linear aggregation task [even if (ii)
enters in our setting], assuming implicitly here that
the models are given in advance and are beyond our control and that the
goal is to combine them appropriately.

Finally, in practice, the noise distribution often departs from the
normal distribution. In particular,
it can exhibit\vadjust{\goodbreak} much heavier tails, and consequently induce highly
non-Gaussian residuals. It is then natural
to ask whether classical estimators such as the ridge regression and
the ordinary least squares estimator
are sensitive to this type of noise, and whether we can design more
robust estimators.

\subsection*{\texorpdfstring{Outline and contributions.}{Outline and
contributions}}

Section \ref{secridge} provides a new analysis of the ridge estimator
and the ordinary least squares estimator, and their variants.\vadjust{\goodbreak}
Theorem \ref{thhfrlam} provides an asymptotic result for the ridge estimator,
while Theorem~\ref{thermom} gives a nonasymptotic risk bound
for the empirical risk minimizer, which is complementary to the theorems
put in the survey section. In particular, the result has the benefit to
hold for the ordinary least squares estimator and for heavy-tailed outputs.
We show quantitatively that the ridge penalty
leads to an implicit reduction of the input space dimension.
Section \ref{seccomputable} shows a nonasymptotic $d/n$ exponential deviation
risk bound under weak moment conditions on the output $Y$ and on the
$d$-dimensional input representation~$\vp(X)$.

The main contribution of this paper is to show through a PAC-Bayesian
analysis on truncated differences of losses that the output
distribution does not need to have bounded
conditional exponential moments in order for the excess risk of appropriate
estimators to concentrate exponentially.
Our results tend to say that truncation leads to more robust algorithms.
Local robustness to contamination is usually invoked
to advocate the removal of outliers,
claiming that estimators should be made insensitive to
small amounts of spurious data.
Our work leads to a different theoretical explanation.
The observed points having unusually large outputs when compared
with the (empirical) variance should be down-weighted in the estimation
of the mean, since they contain less information than noise.
In short, huge outputs should be truncated because of their low
signal-to-noise ratio.

\section{\texorpdfstring{Ridge regression and empirical risk minimization.}{Ridge regression and empirical risk minimization}} \label{secridge}

We recall the definition
\[
\cF= \Biggl\{ f_{\theta} = \sum_{j=1}^d \th_j \vp_j ; (\th
_1,\ldots,\th_d) \in\cC\Biggr\},
\]
where $\cC$ is a closed convex set, not necessarily bounded
(so that $\Theta= \mathbb{R}^d$ is allowed).
In this section we provide exponential deviation inequalities for the
empirical risk
minimizer and the ridge regression estimator on $\cF$ under weak conditions
on the tail of the output distribution.\vadjust{\goodbreak}

The most general theorem which can be obtained from
the route followed in this section is Theorem
1.5 of the supplementary material \cite{AudCat11supp}.
It is expressed in terms of a series of empirical
bounds. The first deduction we can make from this
technical result is of an asymptotic nature.
It is stated under weak hypotheses, taking advantage
of the weak law of large numbers.
\begin{theorem} \label{thhfrlam}
For $\lam\ge0$, let $\frid$ be its associated optimal ridge function
[see~(\ref{eqfrid})].
Let us assume that
%
%
\begin{equation}
\mathbb{E}[ \Vert\vp(X) \Vert^4 ] < + \infty
\end{equation}
and
\begin{equation}
\mathbb{E} \{ \Vert\vp(X) \Vert^2
[\frid(X) - Y ]^2 \} < + \infty.
\end{equation}
Let $\nu_1 > \cdots> \nu_d$ be the eigenvalues of
the Gram matrix $Q=\mathbb{E} [ \vp(X) \vp(X)^T]$, and
let $Q_{\lambda} = Q + \lambda I$ be the ridge regularization of $Q$.
Let us define the \textit{effective ridge dimension}
\[
D = \sum_{i=1}^d \frac{\nu_i}{\nu_i+\lam} \mathbh{1}(\nu_i>0)
= \Tr[ (Q+ \lambda I)^{-1} Q ]
= \mathbb{E} [ \Vert Q_{\lambda}^{-1/2} \vp(X)
\Vert^2 ].
\]
When $\lambda= 0$, $D$ is equal to the rank
of $Q$ and is otherwise smaller.
For any $\varepsilon> 0$, there is $n_{\varepsilon}$,
such that for any $n \geq n_{\varepsilon}$,
with probability at least $1 - \varepsilon$,
\begin{eqnarray*}
&&R\bigl(\hfrlam\bigr)  + \lambda\bigl\Vert\hat{\th}^{(\mathrm{ridge})}
\bigr\Vert^2 \\
&&\qquad \leq\min_{\theta\in\Theta} \{ R(f_{\theta})
+ \lambda\Vert\theta\Vert^2 \} \\
&&\qquad\quad{} + \frac{30 \mathbb{E} \{ \Vert Q_\lam^{-1/2} \vp
(X) \Vert^2
[ \frid(X) - Y ]^2\}}{
\mathbb{E} \{ \Vert Q_\lam^{-1/2} \vp(X) \Vert^2 \}}
\frac{D }{n} \\
&&\qquad\quad{} + 1\mbox{,}000 \sup_{
v \in\mathbb{R}^d}
\frac{\mathbb{E} [ \langle v, \vp(X) \rangle^2 [ \frid
(X) - Y ]^2
]}{\mathbb{E} ( \langle v, \vp(X) \rangle^2 ) +
\lambda\Vert v \Vert^2}
\frac{\log(3\varepsilon^{-1})}{n}\\
&&\qquad \leq\min_{\theta\in\Theta} \{ R(f_{\theta})
+ \lambda\Vert\theta\Vert^2 \} \\
&&\qquad\quad{}
+ \ess\sup\E\{[Y-\frid(X)]^2 | X\} \frac
{30 D + 1\mbox{,}000 \log( 3\varepsilon^{-1} )}{n}.
\end{eqnarray*}
\end{theorem}
\begin{pf}
See Section 1 of the supplementary material \cite{AudCat11supp}.
\end{pf}

This theorem shows that the ordinary least squares estimator
(obtained when $\Theta= \mathbb{R}^d$ and $\lambda= 0$),
as well as the empirical risk minimizer on any closed convex
set, asymptotically reaches a $d/n$ speed of convergence
under very weak hypotheses. It shows also the regularization
effect of the ridge regression. There emerges an \textit{effective
dimension} $D$,
where the ridge penalty has a threshold effect on the eigenvalues
of the Gram matrix.

Let us remark that the second inequality stated in the theorem
provides a simplified bound which makes sense only when
\[
\ess\sup\mathbb{E} \{ [ Y - \frid(X) ]^2 | X
\}
< + \infty
\]
implying that $\Vert\frid- \freg\Vert_{\infty} < + \infty$.
We chose to state the first inequality as well, since it
does not require such a tight relationship between $\frid$ and
$\freg$.\vadjust{\goodbreak}

On the other hand, the weakness of this result is
its asymptotic nature: $n_\varepsilon$ may be arbitrarily large
under such weak hypotheses, and this happens even in the
simplest case of the estimation of the mean of a real-valued
random variable by its empirical mean [which is the case when
$d = 1$ and $\vp(X) \equiv1$].

Let us now give some nonasymptotic rate under stronger
hypotheses and for the empirical risk minimizer (i.e.,
$\lambda= 0$).
\begin{theorem} \label{thermom}
Assume that $\E\{[Y-f^*(X)]^4\}< +\infty$
and
\[
\bcR= \sup_{f\in\Span\{\vp_1,\ldots,\vp_d\}-\{0\}} {\|f\| _\infty
^2}/{\E[f(X)^2]}
< +\infty.
\]
Consider the (unique) empirical risk minimizer $\hferm=f_{\therm}\dvtx x
\mapsto\langle\therm,\break\vp(x)\rangle$ on $\cF$ for which
$\therm\in\Span\{\vp(X_1),\ldots,\vp(X_n)\}$.\footnote{When
$\cF=\Flin$, we have $\therm= \XX^+ \YY$,
with $\XX=(\vp_j(X_i))_{ 1\le i \le n, 1\le j \le d}$, $\YY=[Y_j]_{j=1}^n$
and $\XX^+$ is the Moore--Penrose pseudoinverse of $\XX$.}
For any values of~$\varepsilon$ and $n$ such that
$2/n \le\varepsilon\le1$ and
\[
n > 1280 B^2 \biggl[ 3 B d+ \log(2/\varepsilon) + \frac{16 B^2
{d}^2}{n} \biggr]
\]
with probability at least $1-\varepsilon$,
%
%
\begin{eqnarray}\label{eqdp}
&&R\bigl(\hferm\bigr) - R(f^*) \nonumber\\[-8pt]\\[-8pt]
&&\qquad\le1920 \bcR\sqrt{\E
\{ [Y-f^*(X)]^4 \}} \biggl[ \frac{3 B d + \log
(2\varepsilon
^{-1})}{n} + \biggl(\frac{4 B d}n\biggr)^2 \biggr].
\nonumber
\end{eqnarray}
\end{theorem}
\begin{pf}
See Section 1 of the supplementary material \cite{AudCat11supp}.
\end{pf}

It is quite surprising that the traditional assumption of uniform
boundedness of the conditional exponential moments of the output
can be replaced by a simple moment condition for reasonable confidence
levels (i.e., $\varepsilon\geq2/n$). For highest confidence levels,
things are more tricky since we need to control with high probability a
term of order $[r(f^*)-R(f^*)]d/n$ (see Theorem~1.6).
The cost to pay to get the exponential deviations under only a
fourth-order moment condition on the output is the appearance of the
geometrical quantity $\bcR$ as a multiplicative factor.

To better understand the quantity $\bcR$, let us consider two cases.
First, consider that
the input is uniformly distributed on $\X=[0,1]$, and
that the functions $\vp_1,\ldots,\vp_d$ belong to the Fourier basis.
Then the quantity $\bcR$
behaves like a numerical constant. On the contrary, if we take $\vp
_1,\ldots,\vp_d$ as the first $d$ elements of a wavelet expansion,
the more localized wavelets induce high values of $\bcR$, and $\bcR$
scales like $\sqrt{d}$, meaning that Theorem \ref{thermom}
fails to give a $d/n$-excess risk bound in this case. This limitation
does not appear in Theorem \ref{thhfrlam}.\vadjust{\goodbreak}

To conclude, Theorem \ref{thermom} is limited in at least four ways:
it involves the quantity~$\bcR$, it applies only to
uniformly bounded $\vp(X)$, the output needs to have a fourth moment, and
the confidence level should be as great as $\varepsilon\geq2/n$.
These limitations will be addressed in the next section by
considering a more involved algorithm.

\section{\texorpdfstring{A min--max estimator for robust estimation.}{A min--max estimator for robust
estimation}}
\label{seccomputable}

This section provides an alternative to the empirical risk minimizer
with nonasymptotic exponential risk deviations of order $d/n$ for
any confidence level.
Moreover, we will assume only a second-order moment condition
on the output and cover the case of unbounded inputs, the requirement
on $\varphi(X)$
being only a finite fourth-order moment. On the other hand, we assume here
that the set $\Theta$ of the vectors of coefficients is bounded.
The computability of the proposed estimator and numerical experiments
are discussed at the end of the section.

\subsection{\texorpdfstring{The min--max estimator and its theoretical guarantee.}{The min--max estimator and its theoretical
guarantee}}

Let $\alpha>0$, $\lam\ge0$, and consider the truncation function:
\[
\psi(x) = \cases{
- \log(1 - x + x^2/2 ), &\quad $0 \leq x \leq1$, \cr
\log(2), &\quad $x \geq1$, \cr
- \psi(-x), &\quad $x \leq0$.}
\]
For any $\th,\th'\in\Theta$, introduce
\[
\cD(\th,\th')=n \alpha\lam( \|\th\|^2-\|\th'\|^2 )
+ \sum_{i=1}^n \psi\bigl(\alpha[Y_i-f_\th(X_i)]^2-\alpha
[Y_i-f_{\th'}(X_i)]^2\bigr).
\]
We recall that $\frid=f_{\thrid}$ with
$
\thrid\in\arg\min_{\theta\in\Theta} \{ R(f_{\theta})
+ \lambda\Vert\theta\Vert^2 \}
$, and that the effective ridge dimension is defined as
\[
D = \mathbb{E} [ \Vert Q_{\lambda}^{-1/2} \vp(X)
\Vert^2 ]
= \Tr[ (Q+ \lambda I)^{-1} Q ]
= \sum_{i=1}^d \frac{\nu_i}{\nu_i+\lam} \mathbh{1}(\nu_i>0)
\le d,
\]
where $\nu_1 \geq\cdots\geq\nu_d$ are the eigenvalues of the Gram
matrix $Q=\break\mathbb{E} [ \vp(X) \vp(X)^T]$.
Let us assume in this section that
%
%
\begin{equation} \label{eqas1}
\E\{[Y-\frid(X)]^4\}<+\infty,
\end{equation}
and that for any $j\in\{1,\ldots,d\}$,
%
%
\begin{equation} \label{eqas2}
\E[\vp_j(X)^4]<+\infty.
\end{equation}

Define
%
%
\begin{eqnarray}
\label{eqcs}
\cS& = &\{ f\in\Flin\dvtx\E[ f(X)^2 ] = 1 \},
\\
\sigma& = &\sqrt{\E\{[Y-\frid(X)]^2\}}= \sqrt{R(\frid)}
, \\
\chi& = &\max_{f\in\cS} \sqrt{\E[ f(X)^4 ] },\\
\kappa& = &\frac{ \sqrt{\E\{[\vp(X)^TQ_\lam^{-1}\vp
(X)]^2\}}}
{\E[\vp(X)^TQ_\lam^{-1}\vp(X) ]},\\
\kappa' & = &\frac{\sqrt{\E\{[Y-\frid(X)]^4\}}}{\E\{
[Y-\frid(X)]^2\}}
= \frac{\sqrt{\E\{[Y-\frid(X)]^4\}}}{\sigma^2},\\
\label{eqkapp}
T & = &\max_{\th\in\Theta,\th'\in\Theta}
\sqrt{ \lam\|\th-\th'\|^2 + \E\{ [ f_{\th}(X)-f_{\th'}(X) ]^2
\} }.
\end{eqnarray}
\begin{theorem} \label{th31}
Let us assume that (\ref{eqas1}) and (\ref{eqas2}) hold.
For some numerical constants $c$ and $c'$,
for
\[
n > c \kappa\chi D
\]
by taking
%
%
\begin{equation} \label{eqalpha}
\alpha= \frac{1}{2 \chi[ 2 \sqrt{\kappa'} \sigma
+ \sqrt{\chi} T
]^2} \biggl(1 - \frac{c \kappa\chi D}{n} \biggr)
\end{equation}
for any estimator $f_{\hth}$ satisfying $\hth\in\Theta$ a.s., for
any $\varepsilon>0$ and any
$\lam\ge0$,
with probability at least $1-\varepsilon$, we have
\begin{eqnarray*}
R(f_{\hth}) + \lam\Vert\hth\Vert^2 & \le &
\min_{\theta\in\Theta} \{ R(f_{\theta}) + \lambda\Vert
\theta\Vert^2 \}\\
&&{} + \frac1{n\alpha}\Bigl( \max_{\th_1\in\Theta} \cD(\hth
,\th_1)
- \inf_{\th\in{\Theta}} \max_{\th_1\in{\Theta}} \cD(\th,\th
_1) \Bigr)
+ \frac{c \kappa\kappa' D \sigma^2}{n} \\
&&{}+ 8 \chi\biggl(
\frac{\log(\varepsilon^{-1})}{n} +
\frac{c' \kappa^2 D^2}{n^2} \biggr) \frac{ [ 2 \sqrt{\kappa
'}\sigma
+ \sqrt{\chi} T ]^2}{ 1 - {c\kappa\chi D}/{n} }.
\end{eqnarray*}
\end{theorem}
\begin{pf}
See Section 2 of the supplementary material \cite{AudCat11supp}.
\end{pf}

By choosing an estimator such that
\[
\max_{\th_1\in\Theta} \cD(\hth,\th_1) < \inf_{\th\in{\Theta}}
\max_{\th_1\in{\Theta}} \cD(\th,\th_1) + \sigma^2 \frac{D}n,
\]
Theorem \ref{th31} provides a nonasymptotic bound for the excess
(ridge) risk with a $D/n$ convergence rate and an exponential tail even
when neither the output $Y$ nor the input vector $\vp(X)$ have
exponential moments. This stronger nonasymptotic bound compared to the
bounds of the previous section comes
at the price of replacing the empirical risk
minimizer by a more involved estimator. Section \ref{SECCOMPUT} provides
a way of computing it approximately.

Theorem \ref{th31} requires a fourth-order moment condition on the output.
In fact, one can replace (\ref{eqas1})
by the following\vadjust{\goodbreak} second-order moment condition on the output:
for any $j\in\{1,\ldots,d\}$,
\[
\E\{\vp_j(X)^2[Y-\frid(X)]^2\}<+\infty,
\]
and still obtain a $D/n$ excess risk bound. This comes at the price of
a more lengthy formula, where terms with $\kap'$ become
terms involving the quantities $\max_{f\in\cS} \E\{
f(X)^2[Y-\frid(X)]^2\}$ and
$\E\{\vp(X)^TQ^{-1}\vp(X)[Y-\frid(X)]^2\}$.
(This can be seen by not using Cauchy--Schwarz's inequality in (2.5)
and (2.6) of the supplementary material \cite{AudCat11supp}.)

\subsection{\texorpdfstring{The value of the uncentered kurtosis coefficients $\chi$
and $\kappa$.}{The value of the uncentered kurtosis coefficients chi
and kappa}}
We see that the speed of convergence of the excess risk in
Theorem \ref{th31} (page \pageref{th31}) depends
on three kurtosis-like coefficients, $\chi$, $\kappa$ and $\kappa'$.
The third, $\kappa'$, is concerned with the noise, conceived as
the difference between the observed output $Y$ and its best explanation
$\frid(X)$ according to the ridge criterion.
The aim of this section is to study the order of magnitude of the
two other coefficients $\chi$ and $\kappa$, which are related
to the design distribution,
\[
\chi = \sup
\{ \mathbb{E} ( \langle u, \varphi(X) \rangle^4
)^{1/2} ;
u \in\mathbb{R}^d, \mathbb{E} ( \langle u, \varphi(X) \rangle
^2 ) \leq1 \}
\]
and
\[
\kappa = D^{-1} \mathbb{E} (
\Vert Q_{\lambda}^{-1/2} \varphi(X) \Vert^4 )^{1/2}.
\]
We will review a few typical situations.

\subsubsection{\texorpdfstring{Gaussian design.}{Gaussian design}}
Let us assume first that $\varphi(X)$ is a multivariate centered
Gaussian random variable. In this case, its covariance matrix
coincides with its Gram matrix $Q_{0}$ and can be written as
\[
Q_0 = U^{-1} \Diag(\nu_i, i=1,\ldots, n) U,
\]
where $U$ is an orthogonal
matrix. Using $U$, we can introduce
$W = U Q_{\lambda}^{-1/2} \varphi(X)$.
It is also a Gaussian vector,
with covariance
$\Diag[ \nu_i/(\lambda+ \nu_i), i = 1 ,\ldots, d ]$.
Moreover, since $U$ is orthogonal,
$\Vert W \Vert= \Vert Q_{\lambda}^{-1/2} \varphi(X) \Vert$, and
since $(W_i, W_j)$ are uncorrelated when $i \neq j$, they are
independent, leading to
\begin{eqnarray*}
\mathbb{E} ( \Vert Q_{\lambda}^{-1/2} \varphi(X) \Vert^4
)
&=& \mathbb{E} \Biggl[ \Biggl( \sum_{i=1}^d W_i^2 \Biggr)^2 \Biggr]
\\
&=& \sum_{i=1}^d \mathbb{E} ( W_i^4 ) + 2 \sum_{1 \leq i <
j \leq d}
\mathbb{E} ( W_i^2 ) \mathbb{E} ( W_j^2 )\\
&=& D^2 + 2 D_2,
\end{eqnarray*}
where $ D_2 = \sum_{i=1}^d \frac{\nu_i^2}{( \lambda+ \nu
_i )^2}$.
Thus, in this case,
\[
\kappa= \sqrt{1 + 2 D_2 D^{-2}} \leq\sqrt{1 + \frac{2 \nu
_1}{(\lambda+ \nu_1)D}} \leq\sqrt{3}.\vadjust{\goodbreak}
\]
Moreover, as for any value of $u$,
$\langle u , \varphi(X) \rangle$ is a Gaussian random variable,
$\chi= \sqrt{3}$.

This situation arises in compressed sensing using random projections on
Gaussian vectors.
Specifically, assume that we want to recover a signal $f\in\R^M$ that
we know to be well approximated
by a linear combination of $d$ basis vectors $f_1,\ldots,f_d$. We
measure $n \ll M$ projections of
the signal $f$ on i.i.d. $M$-dimensional standard normal random vectors
$X_1,\ldots,X_n\dvtx Y_i=\langle f,X_i\rangle$, $i=1,\ldots,n$.
Then, recovering the coefficient $\th_1,\ldots,\th_d$ such that
$f=\sum_{j=1}^d \th_j f_j$
is associated to the least squares regression problem,
$Y \approx\sum_{j=1}^d \th_j \varphi_j(X),$
with $\varphi_j(x)=\langle f_j,x\rangle$, and $X$ having a
$M$-dimensional standard normal distribution.

\subsubsection{\texorpdfstring{Independent design.}{Independent design}}

Let us study now the case when almost surely $\varphi_1(X) \equiv1$ and
$\varphi_2(X),\ldots, \varphi_d(X)$ are independent. To compute $\chi$,
we can assume without loss of generality that $\varphi _2(X),\ldots,
\varphi_d(X)$ are centered and of unit variance, since this
renormalization is precisely the linear transformation that turns the
Gram matrix into the identity matrix. Let us introduce
\[
\chi_* = \max_{j=1,\ldots, d} \frac{ \mathbb{E} [ \varphi
_j(X)^4 ]^{1/2}}{\mathbb{E}
[ \varphi_j(X)^2 ]}
\]
with the convention $\frac00=0$.
A computation similar to the one made in the Gaussian case shows that
\[
\kappa\leq\sqrt{1 + (\chi_*^2-1) D_2 D^{-2}} \leq
\sqrt{ 1 + \frac{(\chi_*^2 - 1) \nu_1}{(\lambda+ \nu_1) D}} \leq
\chi_*.
\]
Moreover, for any $u\in\R^d$ such that $\|u\|=1$,
\begin{eqnarray*}
\mathbb{E} ( \langle u, \varphi(X) \rangle^4 ) &=&
\sum_{i=1}^d u_i^4 \mathbb{E} (\varphi_i(X)^4) + 6 \sum_{1\leq i<j
\leq d}
u_i^2 u_j^2 \mathbb{E}[ \varphi_i(X)^2 ] \mathbb{E}
[
\varphi_j(X)^2 ] \\
&&{}+ 4 \sum_{i=2}^d u_1 u_i^3 \mathbb{E}
[ \varphi_i(X)^3 ]
\\
&\leq&\chi_*^2 \sum_{i=1}^d u_i^4
+ 6 \sum_{i < j} u_i^2 u_j^2 + 4 \chi_*^{3/2} \sum_{i=2}^d
\vert u_1 u_i \vert^3
\\
&\leq&\sup_{u \in\R_{+}^d, \Vert u \Vert= 1} ( \chi_*^2-3
) \sum_{i=1}^d u_i^4 + 3 \Biggl( \sum_{i=1}^d
u_i^2 \Biggr)^2 + 4 \chi_*^{3/2} u_1 \sum_{i=2}^d u_i^3
\\
&\leq&\frac{3^{3/2}}{4} \chi_*^{3/2} + \cases{
\chi_*^2 , &\quad $\chi_*^2 \geq3$,\vspace*{2pt}\cr
3 + \dfrac{\chi_*^2 - 3}{d}, &\quad $1 \leq\chi_*^2 < 3$.}
\end{eqnarray*}
Thus, in this case,
\[
\chi\leq\cases{
\chi_* \biggl( 1 + \dfrac{3^{3/2}}{4 \sqrt{\chi_*}} \biggr)^{1/2}, &\quad
$\chi_* \geq\sqrt{3}$,\vspace*{2pt}\cr
\biggl(3 + \dfrac{3^{3/2}}{4} \chi_*^{3/2} + \dfrac{\chi_*^2 - 3}{d}
\biggr)^{1/2}, &\quad $1 \leq\chi_* < \sqrt{3}$.}
\]

If, moreover, the random variables $\varphi_2(X),\ldots,\varphi_d(X)$
are not skewed, in the sense
that $\mathbb{E} [ \varphi_j(X)^3 ] = 0$, $j = 2,\ldots, d$,
then
\[
\cases{
\chi= \chi_*, &\quad $\chi_* \geq\sqrt{3}$,\vspace*{2pt}\cr
\chi\leq\biggl( 3 + \dfrac{\chi_*^2 - 3}{d} \biggr)^{1/2}, &\quad $1
\leq\chi_* <
\sqrt{3}$.}
\]

\subsubsection{\texorpdfstring{Bounded design.}{Bounded design}}

Let us assume now that the distribution of $\varphi(X)$ is almost
surely bounded and nearly orthogonal. These hypotheses are suited to
the study of regression in usual function bases, like the Fourier
basis, wavelet bases, histograms or splines.

More precisely, let us assume that $\mathbb{P}( \Vert\varphi(X)
\Vert\leq B) = 1$ and
that for some positive constant $A$ and any $u \in\mathbb{R}^d$,
\[
\Vert u \Vert\leq A \mathbb{E} [ \langle u, \varphi(X) \rangle
^2 ]^{1/2}.
\]
This appears as some stability property of the partial basis $\varphi_j$
with respect to the $\mathbb{L}_2$-norm, since it can also be written as
\[
\sum_{j=1}^d u_j^2 \leq A^2 \mathbb{E} \Biggl[ \Biggl( \sum
_{j=1}^d u_j
\varphi_j(X) \Biggr)^2 \Biggr],\qquad u \in\mathbb{R}^d.
\]
In terms of eigenvalues, $A^{-2}$ can be taken to be
the lowest eigenvalue $\nu_d$ of the Gram matrix $Q$.
The value of $A$ can also be deduced from a condition saying that
$\varphi_j$ are nearly orthogonal in the sense that
\[
\mathbb{E} [ \varphi_j(X)^2 ] \geq1\quad \mbox{and}
\quad\vert\mathbb{E} [ \varphi_j(X) \varphi_k(X) ]
\vert\leq\frac{1 - A^{-2}}{d-1}.
\]
In this situation, the chain of inequalities
\[
\mathbb{E} [ \langle u, \varphi(X) \rangle^4 ]
\leq\Vert u \Vert^2 B^2 \mathbb{E}
[ \langle u, \varphi(X) \rangle^2 ] \leq A^2 B^2 \mathbb
{E} [
\langle u, \varphi(X) \rangle^2 ]^2
\]
shows that $\chi\leq A B$. On the other hand,
\begin{eqnarray*}
&&\mathbb{E} [ \Vert Q_{\lambda}^{-1/2} \varphi(X) \Vert^4
]\\[1pt]
&&\qquad= \mathbb{E} [ \sup\{ \langle u, \varphi(X) \rangle^4
; u \in\mathbb{R}^d,
\Vert Q_{\lambda}^{1/2} u \Vert\leq1 \} ] \\[1pt]
&&\qquad\leq\mathbb{E} [ \sup\{ \Vert u \Vert^2 B^2 \langle
u, \varphi(X) \rangle^2 ; \Vert Q^{1/2}_{\lambda} u \Vert
\leq1 \} ] \\[1pt]
&&\qquad\leq\mathbb{E} [ \sup\{ ( 1 + \lambda A^2
)^{-1} A^2 B^2
\Vert Q_{\lambda}^{1/2} u \Vert^2 \langle u,
\varphi(X) \rangle^2 ; \Vert Q_{\lambda}^{1/2} u \Vert\leq1
\} ] \\[1pt]
&&\qquad\leq\frac{A^2B^2}{1 + \lambda A^2} \mathbb{E} [ \Vert
Q_{\lambda}^{-1/2}
\varphi(X) \Vert^2 ] = \frac{A^2 B^2 D}{1 + \lambda A^2}
\end{eqnarray*}
showing that $\kappa\leq\frac{A B}{\sqrt{(1 + \lambda A^2) D}}$.

For example, if $X$ is the uniform random variable on the unit interval
and~$\varphi_j$, $j=1,\ldots, d$, are any functions from the Fourier
basis [meaning that they are of the form $\sqrt{2} \cos(2 k \pi X)$
or $\sqrt{2} \sin( 2 k \pi X)$],
then $A=1$ (because they form an orthogonal
system) and $B \leq\sqrt{2d}$.

A localized basis like the evenly spaced histogram
basis of the unit interval
\[
\varphi_j(x) =
\sqrt{d} \mathbh{1}\bigl(x \in[(j-1)/d, j/d [
\bigr),\qquad j=1,\ldots,
d,
\]
will also be such that $A = 1$ and $B = \sqrt{d}$.
Similar computations could be made for other local bases, like
wavelet bases.

Note that when $\chi$ is of order $\sqrt{d}$, and $\kappa$
and $\kappa'$ of order $1$, Theo-\break rem~\ref{th31}
means that the excess risk of the min--max truncated estimator~%
$\hf$ is upper bounded by $Cd/n$
provided that $n\ge C d^2$ for a large enough constant~$C$.\looseness=1

\subsubsection{\texorpdfstring{Adaptive design planning.}{Adaptive design
planning}}

Let us discuss the case when $X$ is some observed random variable whose
distribution is only approximately known. Namely, let us assume that
$(\varphi_j)_{j=1}^d$ is some basis of functions in $\mathbb{L}_2 [
\wt{\mathbb{P}} ]$ with some known coefficient $\wt{\chi}$, where
$\wt{\mathbb{P}}$ is an approximation of the true distribution of $X$
in the sense that the density of the true distribution $\P$ of $X$ with
respect to the distribution $\wt{\mathbb{P}}$ is in the range
$(\eta^{-1}, \eta)$. In this situation, the coefficient $\chi$
satisfies the inequality $\chi\leq\eta^{3/2} \wt{\chi}$. Indeed,
\begin{eqnarray*}
\mathbb{E}_{X\sim\P} [ \langle u, \varphi(X) \rangle^4 ]
& \leq &\eta{\mathbb{E}}_{X\sim\wt{\P}} [ \langle u, \varphi
(X) \rangle^4 ] \\[1pt]
& \leq &\eta\wt{\chi}^2 {\mathbb{E}}_{X\sim\wt{\P}} [
\langle u, \varphi(X) \rangle^2 ]^2\\[1pt]
&\leq&\eta^3 \wt{\chi}^2 \mathbb{E}_{X\sim{\P}} [ \langle u,
\varphi(X) \rangle^2 ]^2.
\end{eqnarray*}
In the same way, $\kappa\leq\eta^{7/2} \wt{\kappa}$. Indeed,
\begin{eqnarray*}
&&
\mathbb{E} [ \sup\{ \langle u, \varphi(X) \rangle^4 ;
\mathbb{E} ( \langle u, \varphi(X) \rangle^2 ) \leq1
\} ]
\\
&&\qquad\leq\eta\wt{\mathbb{E}} [ \sup\{ \langle u, \varphi
(X) \rangle^4
; \wt{\mathbb{E}} ( \langle u, \varphi(X) \rangle^2 )
\leq\eta\} ] \\
&&\qquad\leq\eta^3 \wt{\mathbb{E}} [ \sup\{ \langle u,
\varphi(X) \rangle^4 ; \wt{\mathbb{E}} ( \langle u,
\varphi(X) \rangle^2
) \leq1 \} ] \\
&&\qquad\leq\eta^3 \wt{\kappa}^2 \wt{\mathbb{E}} [ \sup\{
\langle u,
\varphi(X) \rangle^2 ; \wt{\mathbb{E}} ( \langle u, \varphi(X)
\rangle^2 ) \leq1 \} ]^2 \\
&&\qquad\leq\eta^7 \wt{\kappa}^2 \mathbb{E} [ \sup
\{ \langle u, \varphi(X) \rangle^2 ; \mathbb{E}
( \langle u, \varphi(X) \rangle^2 ) \leq1 \} ]^2.
\end{eqnarray*}

Let us conclude this section with some scenario
for the case when $X$ is a~real-valued random variable. Let us consider the
distribution function of~$\wt{\mathbb{P}}$,
\[
\wt{F}(x) = \wt{\mathbb{P}} ( X \leq x).\vadjust{\goodbreak}
\]
Then, if $\wt{\mathbb{P}}$ has no atoms,
the distribution of $\wt{F}(X)$ would be uniform on $(0,1)$ if
$X$ were distributed according to $\wt{\mathbb{P}}$. In other words,
$\wt{\mathbb{P}} \circ\wt{F}^{-1} = \mathbb{U}$, the uniform distribution
on the unit interval. Starting from some
suitable partial basis $(\varphi_j)_{j=1}^d$ of $\mathbb{L}_2 [
(0,1), \mathbb{U} ]$ like the
ones discussed above, we can build a basis for our problem as
\[
\wt{\varphi}_j(X) = \varphi_j [ \wt{F}(X) ].
\]
Moreover, if $\mathbb{P}$ is absolutely continuous with respect to
$\wt{\mathbb{P}}$ with
density $g$, then $\mathbb{P} \circ\wt{F}^{-1}$ is absolutely
continuous with
respect to $\wt{\mathbb{P}} \circ\wt{F}^{-1} = \mathbb{U}$, with
density $g \circ\wt{F}^{-1}$,
and, of course, the fact that $g$ takes values in $(\eta^{-1}, \eta)$ implies
the same property for $g \circ\wt{F}^{-1}$. Thus, if $\wt{\chi}$ and
$\wt{\kappa}$ are the coefficients corresponding to $\varphi_j(U)$
when $U$
is the uniform random variable on the unit interval, then the true
coefficient $\chi$ [corresponding to
$\wt{\varphi}_j(X)$] will be such that $\chi\leq\eta^{3/2} \wt
{\chi}$
and $\kappa\leq\eta^{7/2} \wt{\kappa}$.

\subsection{\texorpdfstring{Computation of the estimator.}{Computation of the estimator}} \label{SECCOMPUT}

For ease of description of the algorithm, we will write $X$ for $\vp
(X)$, which is equivalent to considering without loss of generality
that the input space is $\R^d$ and that the functions $\vp_1,\ldots
,\vp_d$ are the coordinate functions.
Therefore, the function $f_\th$ maps an input $x$ to $\langle\theta,
x \rangle$.
Let us introduce
\[
\bL_i(\th) = \alpha( \langle\theta, X_i \rangle- Y_i )^2.
\]
For any subset of indices $I \subset\{1,\ldots, n\}$, let us define
\[
r_I(\theta) = \lam\|\th\|^2 + \frac{1}{\alpha\vert I \vert} \sum
_{i \in I} \bL_i(\th).
\]

We suggest the following heuristics to compute an approximation of
\[
\mathop{\arg\min}_{\theta\in\Theta}
\sup_{\th' \in\Theta} \cD(\theta, \th')\mbox{:}
\]
\begin{itemize}
\item Start from $I_1 = \{1,\ldots, n\}$
with the ordinary least squares estimate
\[
\wh{\theta}_1 = \mathop{\arg\min}_{\R^d} r_{I_1}.
\]
\item At step number $k$, compute
\[
\wh{Q}_k = \frac{1}{\vert I_k \vert} \sum_{i \in I_k} X_i X_i^T.
\]
\item Consider the sets
\[
J_{k,1}(\eta) = \bigl\{ i \in I_k \dvtx\bL_i(\wh{\theta}_k)
X_i^T \wh{Q}_{k}^{-1} X_i \bigl(
1 + \sqrt{ 1 + [\bL_i(\wh{\theta}_k)]^{-1}} \bigr)^2
< \eta\bigr\},
\]
where $\wh{Q}_{k}^{-1}$ is the (pseudo-)inverse of the matrix $\wh{Q}_{k}$.
\item Let us define
\begin{eqnarray*}
\theta_{k,1}(\eta) & = &\mathop{\arg\min}_{\R^d} r_{J_{k,1}(\eta)}, \\
J_{k,2}(\eta) & = & \{ i \in I_k \dvtx\vert\bL_i(\theta
_{k,1}(\eta))
- \bL_i( \wh{\theta}_k ) \vert\le1
\}, \\
\theta_{k,2} (\eta) & = & \mathop{\arg\min}_{\R^d} r_{J_{k,2}(\eta)}, \\
(\eta_{k}, \ell_k) & = & \mathop{\arg\min}_{\eta\in\mathbb{R}_+, \ell\in\{
1, 2\}} \max_{j=1,\ldots, k}
\cD( \theta_{k,\ell}(\eta), \wh{\theta}_j ), \\
I_{k+1} & = & J_{k, \ell_k}(\eta_k), \\
\wh{\theta}_{k+1} & = & \theta_{k, \ell_k}(\eta_k).
\end{eqnarray*}
\item Stop when
\[
\max_{j=1,\ldots, k} \cD(\wh{\theta}_{k+1}, \wh{\theta}_j) \geq0,
\]
and set $\wh{\theta}=\wh{\theta}_k$ as the final estimator of
$\thrid$.

\end{itemize}
Note that there will be at most $n$ steps, since $I_{k+1} \varsubsetneq
I_k$ and in practice much less in this
iterative scheme.
Let us give some justification for this proposal. Let us notice first that
\begin{eqnarray*}
&&\cD(\theta+ h, \theta) \\
&&\qquad= n \alpha\lam(\|\th+ h\|^2-\|\th\|^2)\\
&&\qquad\quad{}+ \sum_{i=1}^n \psi\bigl(
\alpha[ 2 \langle h, X_i \rangle(\langle\theta, X_i
\rangle- Y_i) +
\langle h, X_i \rangle^2 ] \bigr).
\end{eqnarray*}
Hopefully, $\thrid= \arg\min_{\th\in\R^d} (R(f_\th)+\lam\|
\th\|^2)$
is in some small neighborhood
of $\wh{\theta}_k$ already, according to the distance defined by $
Q \simeq\wh{Q}_k$. So we\vspace*{2pt} may try to look for improvements of $\wh
{\theta}_k$
by exploring neighborhoods\vspace*{1pt} of $\wh{\theta}_k$ of increasing
sizes with respect to some approximation of the relevant
norm $\Vert\theta\Vert_Q^2 = \mathbb{E} [ \langle\theta, X
\rangle^2 ]$.

Since the truncation function $\psi$ is constant on $(-\infty,-1]$
and $[1,+\infty)$,
the map $\theta\mapsto\cD(\theta, \wh{\theta}_k)$ induces a
decomposition of
the parameter space into cells corresponding to different sets $I$ of
examples. Indeed, such a set $I$ is associated
to the set $\mathcal{C}_I$ of $\th$ such that $\bL_i(\th)-\bL
_i(\wh{\theta}_k)< 1$ if and only if $i\in I$.
Although this may not be the case, we will do as if
the map $\theta\mapsto\cD(\theta, \wh{\theta}_k)$ restricted to
the cell $\mathcal{C}_I$ reached its minimum at some interior point
of $\mathcal{C}_I$, and approximates this minimizer by the minimizer
of $r_I$.

The idea is to remove first the examples which will become
inactive in the closest cells to the current estimate $\wh{\theta}_k$.
The cells for which the contribution of example number $i$
is constant are delimited by at most four parallel hyperplanes.

It is easy to see that the square of the inverse of the distance of
$\wh{\theta}_k$ to the closest of these hyperplanes is
equal to
\[
\frac1\alpha X_i^T \wh{Q}_{k}^{-1} X_i \bL_i(\wh{\theta}_k) \Biggl(
1 + \sqrt{ 1 + \frac{1}{\bL_i(\wh{\theta}_k)}} \Biggr)^2.
\]
Indeed, this distance is the infimum of $\Vert\wh{Q}_k^{1/2} h \Vert
$, where $h$ is a solution of
\[
\langle h, X_i \rangle^2 + 2 \langle h, X_i \rangle
(\langle\wh{\theta}_k, X_i \rangle- Y_i) = \frac
{1}{\alpha}.
\]
It is computed by considering $h$ of the form $h = \xi\Vert\wh
{Q}_k^{-1/2} X_i \Vert^{-1}
\wh{Q}_k^{-1} X_i$ and solving an equation of order two in $\xi$.

This explains the proposed choice of $J_{k,1}(\eta)$. Then a first
estimate $\theta_{k,1}(\eta)$
is computed on the basis of this reduced sample, and the sample is
readjusted to $J_{k,2}(\eta)$
by checking which constraints are really activated in the computation of
$\cD(\theta_{k,1}(\eta), \wh{\theta}_k)$. The estimated parameter
is then readjusted, taking
into account the readjusted sample (this could as a variant be iterated
more than once).
Now that we have some new candidates $\theta_{k,\ell}(\eta)$, we
check the minimax
property against them to elect $I_{k+1}$ and $\wh{\theta}_{k+1}$.
Since we did not check
the minimax property against the whole parameter set $\Theta=\R^d$,
we have no theoretical
warranty for this simplified algorithm. Nonetheless, similar
computations to what we
did could prove that we are close to solving $\min_{j=1,\ldots, k}
R(f_{\wh{\theta}_j})$,
since we\vspace*{-1pt} checked the minimax property on the reduced parameter set $\{
\wh{\theta}_j,
j=1,\ldots, k \}$. Thus, the proposed heuristics are capable of
improving on the performance
of the ordinary least squares estimator, while being guaranteed not to
degrade its
performance significantly.

\subsection{\texorpdfstring{Synthetic experiments.}{Synthetic
experiments}}

In Section \ref{secnoise}, we detail the different kinds of noises we
work with.
Then, Sections \ref{secexpind}, \ref{sechcc} and \ref{sects} describe
the three types of functional relationships between the input, the
output and the noise involved
in our experiments. A motivation for choosing these input--output
distributions was the ability to compute
exactly the excess risk, and thus to compare easily estimators.
Section \ref{secexper} provides details about the implementation, its
computational efficiency
and the main conclusions of the numerical experiments. Figures and
tables are postponed to the \hyperref[app]{Appendix}.

\subsubsection{\texorpdfstring{Noise distributions.}{Noise distributions}} \label{secnoise}

In our experiments, we consider different ty\-pes of noise that are
centered and with unit variance:
\begin{itemize}
\item the standard Gaussian noise, $W \sim\N(0,1)$,
\item a heavy-tailed noise defined by
$W=\sign(V)/|V|^{1/q}$,
with $V \sim\N(0,1)$, a~standard Gaussian random variable and
$q=2.01$ (the
real number $q$ is taken strictly larger than $2$ as for $q=2$, the
random variable $W$ would not admit
a finite second moment).
\item an asymmetric heavy-tailed noise defined by
\[
W=\cases{
|V|^{-1/q}, &\quad if $V>0$, \vspace*{2pt}\cr
-\dfrac{q}{q-1}, &\quad otherwise,}
\]
with $q=2.01$
with $V \sim\N(0,1)$ a standard Gaussian random variable.
\item a mixture of a Dirac random variable with a low-variance Gaussian
random variable defined by,
with probability $p$, $W=\sqrt{{(1-\rho)}/{p}}$, and with
probability $1-p$, $W$ is drawn from
\[
\N\biggl(-\frac{\sqrt{p(1-\rho)}}{1-p},\frac{\rho}{1-p}-\frac
{p(1-\rho)}{(1-p)^2}\biggr).
\]
The parameter $\rho\in[p,1]$ characterizes the
part of the variance of $W$ explained by the Gaussian part of the mixture.
Note that this noise admits exponential moments, but for $n$ of order
$1/p$, the
Dirac part of the mixture generates low signal-to-noise points.
\end{itemize}

\subsubsection{\texorpdfstring{Independent normalized covariates [$\operatorname{INC}(n,d)$].}{Independent normalized covariates [INC(n,d)]}} \label
{secexpind}
In INC$(n,d)$, we consider $\varphi(X)=X$, and the input--output pair
is such that
\[
Y=\langle\th^* , X \rangle+ \sigma W,
\]
where the components of $X$ are independent standard normal distributions,
$\th^*=(10,\ldots,10)^T \in\R^d$ and $\sigma=10$.

\subsubsection{\texorpdfstring{Highly correlated covariates [$\operatorname{HCC}(n,d)$].}{Highly correlated covariates [HCC(n,d)]}} \label{sechcc}
In $\operatorname{HCC}(n,d)$, we consider $\varphi(X)=X$, and the input--output pair
is such that
\[
Y=\langle\th^* , X \rangle+ \sigma W,
\]
where $X$ is a multivariate centered normal Gaussian with
covariance matrix~$Q$ obtained by drawing a $(d,d)$-matrix $A$ of
uniform random variables
in $[0,1]$\vadjust{\goodbreak} and by computing $Q=A A^T$,
$\th^*=(10,\ldots,10)^T \in\R^d$ and $\sigma=10$.
So the only difference with the setting of Section \ref{secexpind} is
the correlation between the covariates.

\subsubsection{\texorpdfstring{Trigonometric series [$\operatorname{TS}(n,d)$].}{Trigonometric series [TS(n,d)]}} \label{sects}
Let $X$ be a uniform random variable on $[0,1]$.
Let $d$ be an even number.
In TS$(n,d)$, we consider
\[
\varphi(X)= (\cos(2 \pi X),\ldots,\cos(d \pi X),\sin(2 \pi
X),\ldots,\sin(d \pi X))^T,
\]
and the input--output pair is such that
\[
Y=20 X^2-10 X-\tfrac53+\sigma W
\]
with $\sigma=10$.
One can check that this implies
\[
\theta^*=\biggl(\frac{20}{\pi^2},\ldots,\frac{20}{\pi^2 ({d}/2)^2},
-\frac{10}{\pi},\ldots,-\frac{10}{\pi({d}/2)}\biggr)^T
\in\R^d.\vadjust{\goodbreak}
\]

\subsubsection{\texorpdfstring{Experiments.}{Experiments}} \label{secexper}

\paragraph*{\texorpdfstring{Choice of the parameters and implementation details.}{Choice of the parameters and implementation
details}}
The min--max truncated algorithm has two parameters $\alpha$ and
$\lambda$.
In the subsequent experiments, we set the ridge parameter $\lam$ to
the natural default choice for it:
$\lam=0$. For the truncation parameter $\alpha$, according to our
analysis [see (\ref{eqalpha})],
it roughly should be of order $1/\sigma^2$ up to kurtosis coefficients.
By using the ordinary least squares estimator, we roughly estimate this value,
and test values of $\alpha$ in a geometric grid (of $8$ points) around
it (with ratio $3$).
Cross-validation can be used to select the final~$\alpha$.
Nevertheless, it is computationally expensive and is significantly
outperformed in
our experiments by the following simple procedure:
start with the smallest $\alpha$ in the geometric grid and increase it
as long as
$\hat{\th}=\th_1$, that is, as long as we stop at the end of the
first iteration
and output the empirical risk minimizer.

To compute $\theta_{k,1}(\eta)$ or $\theta_{k,2}(\eta)$,
one needs to determine a least squares estimate (for a modified sample).
To reduce the computational burden, we do not want to test all possible
values of
$\eta$ (note that there are at most~$n$ values leading to different estimates).
Our experiments show that testing only three levels of $\eta$ is sufficient.
Precisely, we sort the quantity
\[
\bL_i(\wh{\theta}_k)
X_i^T \wh{Q}_{k}^{-1} X_i \bigl(
1 + \sqrt{ 1 + [\bL_i(\wh{\theta}_k)]^{-1}} \bigr)^2
\]
by decreasing order
and consider $\eta$ being the first, $5$th and $25$th value of the
ordered list.
Overall, in our experiments, the computational complexity is
approximately fifty times
larger than the one of computing the ordinary least squares estimator.

\begin{figure}

\includegraphics{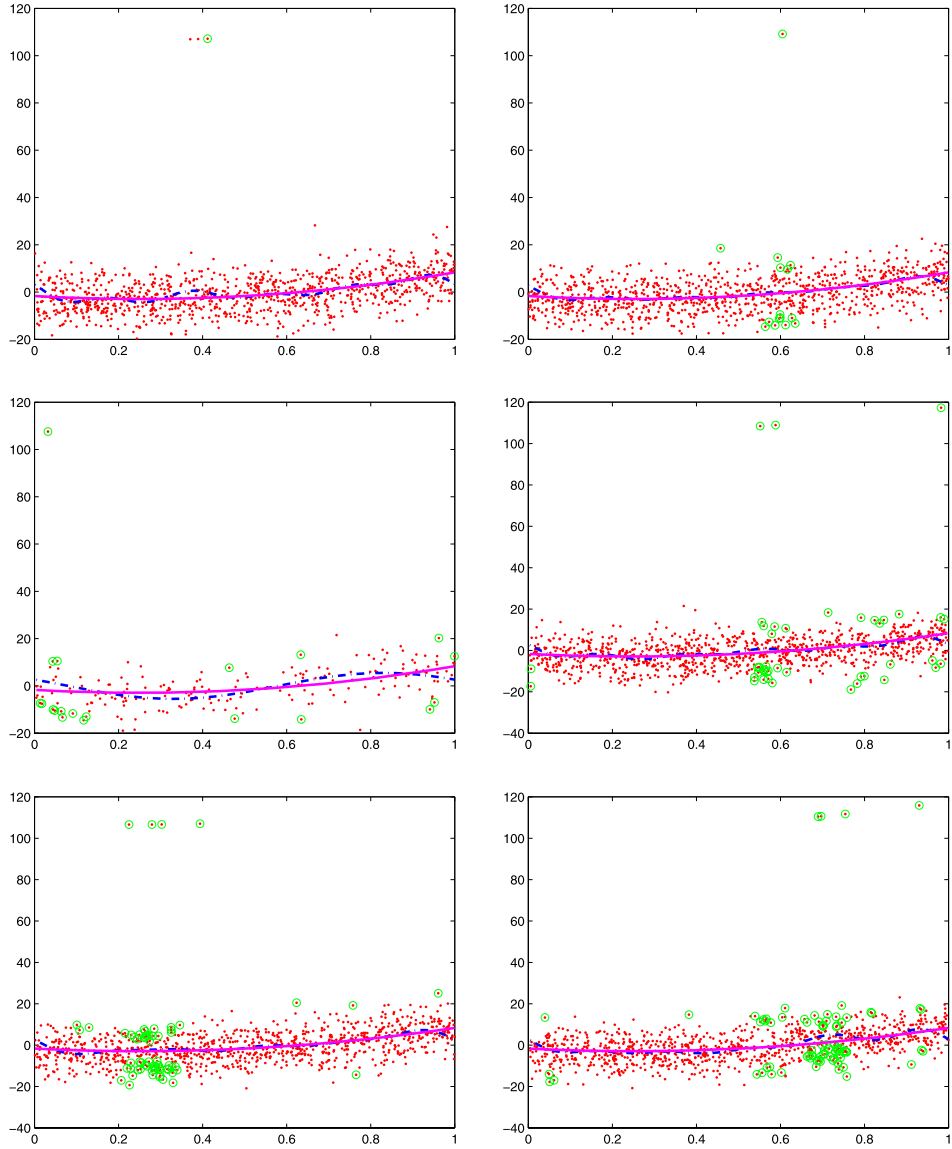}

\caption{Circled points are the points of the training set generated
several times from TS$(1\mbox{,}000,10)$ (with the mixture noise with $p=0.005$
and $\rho=0.4$) that are not taken into account in the min--max
truncated estimator (to the extent that the estimator would not change
by removing simultaneously all these points).
The min--max truncated estimator $x\mapsto\hf(x)$ appears in dash-dot
line, while $x\mapsto\E(Y|X=x)$ is in solid line. In these six simulations,
it outperforms the ordinary least squares estimator.} \label{fig1}
\end{figure}

\begin{figure}

\includegraphics{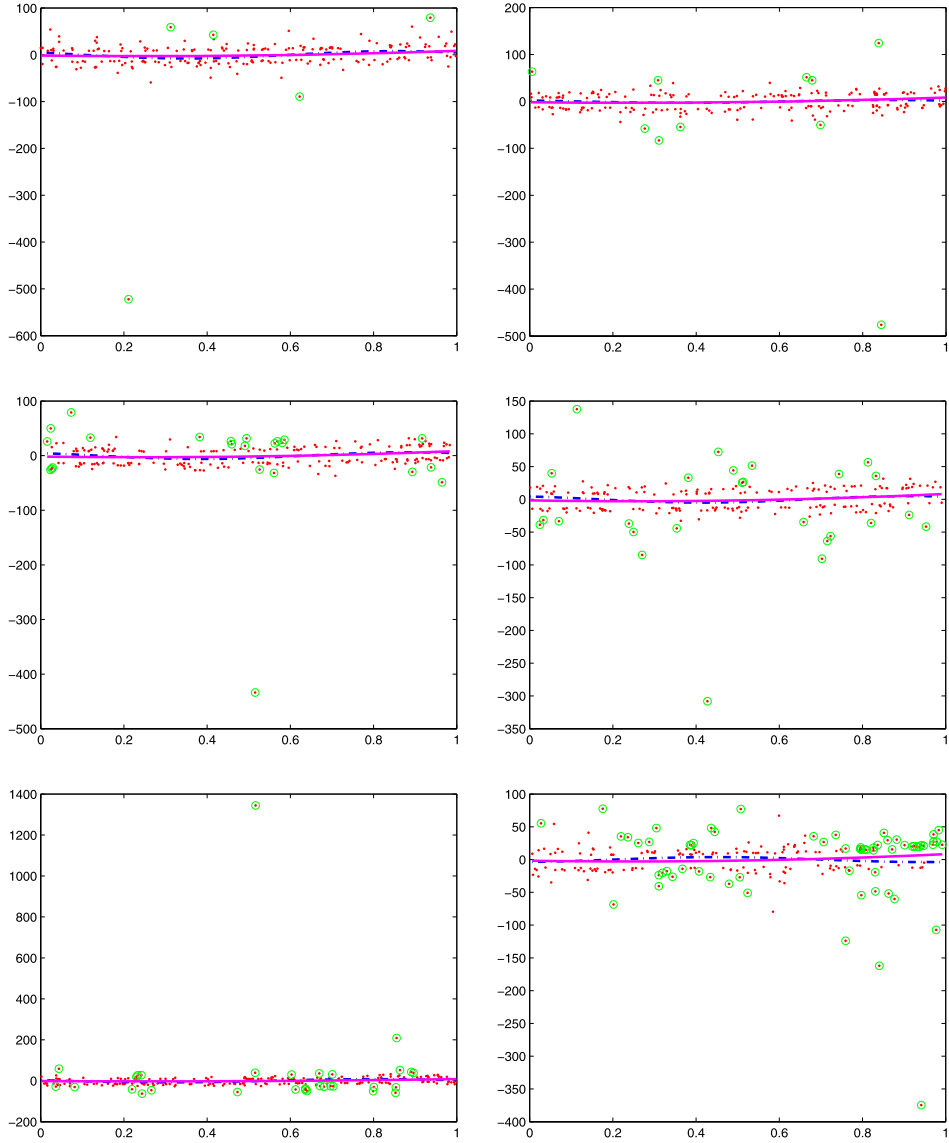}

\caption{Circled points are the points of the training set generated
several times from TS$(200,2)$ (with the heavy-tailed noise) that are
not taken into account in the min--max truncated estimator (to the
extent that the estimator would not change by removing these points).
The min--max truncated estimator $x\mapsto\hf(x)$ appears in dash-dot
line, while $x\mapsto\E(Y|X=x)$ is in solid line. In these six
simulations, it outperforms the ordinary least squares estimator. Note
that in the last figure, it does not consider $64$ points among the
$200$ training points.} \label{fig2}
\end{figure}

\paragraph*{\texorpdfstring{Results.}{Results}}
The tables and figures have been gathered in the \hyperref[app]{Appendix}.
Tables~\ref{tabb01} and \ref{tabb04} give the results for the
mixture noise.
Tables \ref{taba201}, \ref{taba-201} and \ref{taba0} provide the
results for the heavy-tailed noise
and the standard Gaussian noise. Each line of the tables has been
obtained after 1,000 generations of the training set.\vadjust{\goodbreak}
These results show that the min--max truncated estimator is often equal
to the ordinary least squares estimator $\hfols$,
while it ensures impressive consistent improvements when it differs
from $\hfols$.
In this latter case, the number of points that are not considered in
$\hf$, that is, the number of points
with low signal-to-noise ratio, varies a lot from $1$ to $150$ and is
often of order $30$. Note that not only the points
that we expect to be considered as outliers (i.e., very large output
points) are erased, and
that these points seem to be taken out by local groups: see Figures \ref{fig1}
and \ref{fig2}
in which the erased points are marked by surrounding circles.

Besides, the heavier the noise tail is (and also the larger the
variance of the noise is), the more often the truncation modifies the
initial ordinary least squares estimator, and the more improvements we
get from the min--max truncated estimator, which also becomes much more
robust than the ordinary least squares estimator (see the confidence
intervals in the tables).

Finally, we have also tested more traditional methods in robust
regression, namely, the M-estimators with Huber's loss, $L_1$-loss and
Tukey's bisquare influence function, and also the least trimmed squares
estimator, the S-es\-timator and the MM-estimator (see \cite
{Rou84,Yoh87} and the references within). These methods rely on
diminishing the influence of points having ``unreasonably'' large
residuals. They
were developed to handle training sets containing true outliers, that
is, points $(X,Y)$ not generated by the distribution $P$. This is not
the case in our estimation framework.
By overweighting points having reasonably small residuals, these
methods are often biased even in settings where
the noise is symmetric and the regression function $\freg\dvtx x\mapsto
\E
[ Y | X=x]$ belongs to $\Flin$ (i.e., $\freg=f^*_{\mathrm
{lin}}$), and also even when there is no noise (but $\freg\notin
f^*_{\mathrm{lin}}$).

The worst results were obtained by the $L_1$-loss, since estimating the
(conditional) median is here really different from estimating the
(conditional) mean. The MM-estimator and the M-estimators with Huber's
loss and Tukey's bisquare influence function give good results as long
as the signal-to-noise ratio is low. When the signal-to-noise ratio is high,
a lack of consistency drastically appears in part of our simulations,
showing that these methods are thus not suited for our estimation framework.

The S-estimator is almost consistently improving on the ordinary least
squares estimator (in our simulations). However, when the
signal-to-noise ratio is low (i.e., in the setting of the
aforementioned simulations with \mbox{$\sigma=10$}), the improvements are
much less significant than the ones of the min--max truncated estimator.

\section{\texorpdfstring{Main ideas of the proofs.}{Main ideas of the
proofs}}

The goal of this section is to explain the key ingredients appearing in
the proofs which both allow to obtain subexponential tails
for the excess risk under a nonexponential moment assumption and get
rid of the logarithmic factor in the excess risk bound.

\subsection{\texorpdfstring{Subexponential tails under a nonexponential moment
assumption via truncation.}{Subexponential tails under a nonexponential moment
assumption via truncation}}

Let us start with the idea allowing us to prove exponential
inequalities under just a moment assumption (instead of the traditional
exponential moment assumption).
To understand it, we can consider the (apparently) simplistic
$1$-dimensional situation in which we have $\Theta=\R$ and the
marginal distribution of $\varphi_1(X)$ is the Dirac distribution at $1$.
In this case, the risk of the prediction function $f_\th$ is $R(f_\th)=
\E[ (Y-\th)^2 ]=\E[(Y-\E Y)^2 ] + (\E Y -\th)^2$,
so that the least squares regression problem boils down to the
estimation of the mean of the output variable.
If we only assume that $Y$ admits a finite second moment, say, $\E
( Y^2
) \leq1$, it is not clear whether for any $\varepsilon>0$, it is
possible to find $\hth$ such that,
with probability at least $1-2\varepsilon$,
%
%
\begin{equation} \label{eqtar1}
R(f_{\hth})-R(f^*) = \bigl(\E(Y) -\hth\bigr)^2 \le c \frac
{\logeps}{n}
\end{equation}
for some numerical constant $c$.
Indeed, from Chebyshev's inequality, the trivial choice $\hth=\frac{1}{n}
\sum_{i=1}^n Y_i$ just satisfies, with probability at least
$1-2\varepsilon$,
\[
R(f_{\hth})-R(f^*) \le\frac1{n\varepsilon},
\]
which is far from the objective (\ref{eqtar1}) for small confidence
levels [consider $\varepsilon=\exp(-\sqrt{n})$, e.g.].
The key idea is thus to average (soft) \textit{truncated} values of the
outputs. This is performed by
taking
\[
\hth= \frac1{n\lam}\sum_{i=1}^n \log\biggl(1+\lam Y_i+\frac{\lam
^2Y_i^2}{2}\biggr)
\]
with $\lam=\sqrt{\frac{2 \logeps}n}$.
Since we have
\begin{eqnarray*}
\log\E\exp(n\lam\hth) &=& n \log\biggl( 1+ \lam\E(Y) + \frac
{\lam^2}2 \E( Y^2 ) \biggr) \\
&\le& n\lam\E( Y )+ n \frac{\lam^2}2,
\end{eqnarray*}
the exponential Chebyshev's inequality guarantees that with probability
at least $1-\varepsilon$, we have
$
n\lam( \hth-\E(Y) ) \le n \lam^2/2 + \logeps
$,
hence,
\[
\hth-\E(Y)\le\sqrt\frac{2\logeps}{n}.
\]
Replacing $Y$ by $-Y$ in the previous argument, we obtain that,
with probability at least $1-\varepsilon$, we have
\[
n\lam\Biggl\{ \E(Y)+\frac1{n\lam}\sum_{i=1}^n \log\biggl(1-\lam
Y_i+\frac{\lam^2Y_i^2}{2}\biggr) \Biggr\} \le n \frac{\lam^2}2 +
\logeps.
\]
Since $-\log(1+x+x^2/2) \le\log(1-x+x^2/2)$, this implies
$\E(Y)-\break\hth\le\sqrt\frac{2\logeps}{n}$.
The two previous inequalities imply inequality (\ref{eqtar1}) (for
$c={2}$), showing that subexponential tails are achievable even when we
only assume that the random variable admits a finite second moment (see
\cite{Cat09} for more details
on the robust estimation of the mean of a random variable).

\subsection{\texorpdfstring{Localized PAC-Bayesian inequalities to eliminate a
logarithm factor.}{Localized PAC-Bayesian inequalities to eliminate a
logarithm factor}}
Let us first recall that the Kullback--Leibler divergence between
distributions~$\rho$ and $\mu$ defined on $\cF$ is
%
%
\begin{equation} \label{eqkl}
K(\rho,\mu) \eqdef\cases{\displaystyle
\undc{\E}{f\sim\rho} \log\biggl[ \frac{d \rho}{d \mu}(f)
\biggr],
&\quad if $\rho\ll\mu$,\vspace*{2pt}\cr
+ \infty, &\quad otherwise,}
\end{equation}
where $\frac{d \rho}{d \mu}$ denotes as usual
the density of $\rho$ w.r.t. $\mu$.
For any real-valued (measurable) function $h$ defined on $\cF$ such
that\vadjust{\goodbreak}
$\int\exp[h(f)] \pi(df)<+\infty$, we define the distribution
$\pi_{h}$ on $\cF$ by its density:
%
%
\begin{equation} \label{eqpih}
\frac{d \pi_{h}}{d \pi}(f) =
\frac{\exp[ h(f)]}{\int\exp[h(f')] \pi(df')}.
\end{equation}

The analysis of statistical inference generally relies on upper
bounding the supremum of an empirical process $\chi$ indexed by the
functions in a model $\cF$.
Concentration inequalities appear as a central tool to obtain these bounds.
An alternative approach, called the PAC-Bayesian one, consists in using
the entropic equality
%
%
\begin{equation} \label{eqiexp}\quad
\E\exp\biggl( \sup_{\rho\in\M} \biggl\{ \int\rho(df) \chi(f) -
K(\rho,\pi') \biggr\} \biggr)= \int\pi'(df) \E\exp( \chi(f)
),
\end{equation}
where $\M$ is the set of probability distributions on $\cF$.

Let $\chr\dvtx\cF\ra\R$ be an observable process such that, for any
$f\in\cF$, we have
\[
\E\exp(\chi(f)) \le1
\]
for $\chi(f) =\lam[ R(f) - \chr(f)]$ and some $\lam>0$. Then (\ref
{eqiexp}) leads to, for any $\varepsilon>0$, with probability at least
$1-\varepsilon$,
for any distribution $\rho$ on $\cF$, we have
%
%
\begin{equation} \label{eqipac}
\int\rho(df) R(f) \le\int\rho(df) \chr(f) + \frac{K(\rho,\pi
')+\logeps}{\lam}.
\end{equation}
The left-hand side quantity represents the expected risk with respect
to the distribution $\rho$. To get the smallest upper bound on this quantity,
a~natural choice of the (posterior) distribution $\rho$ is obtained by
minimizing the right-hand side, that is, by taking $\rho=\pi'_{-\lam
\chr}$
[with the notation introduced in~(\ref{eqpih})]. This distribution
concentrates on functions $f\in\cF$ for which $\chr(f)$ is small.
Without prior knowledge, one may want to choose a prior distribution
$\pi'=\pi$ which is rather ``flat'' (e.g., the one induced by the
Lebesgue measure in the case of
a model $\cF$ defined by a bounded parameter set in some Euclidean space).
Consequently, the Kullback--Leibler divergence $K(\rho,\pi')$, which
should be seen as the complexity term, might be excessively large.

To overcome the lack of prior information and the resulting high
complexity term, one can alternatively use a more ``localized'' prior
distribution.
Here we use Gaussian distributions centered at the function of interest
(e.g., the function $f^*$), and with covariance matrix proportional
to the inverse of the Gram matrix $Q$. The idea of using PAC-Bayesian
inequalities with Gaussian prior and posterior distributions
goes back to Langford and Shawe-Taylor~\cite{LanSha02} in the context
of linear classification.

The detailed proofs of Theorems \ref{thhfrlam}, \ref{thermom} and \ref{th31}
can be found in the supplementary material \cite{AudCat11supp}.

\begin{appendix}\label{app}

\section*{\texorpdfstring{Appendix: Experimental results for the min--max truncated estimator (Section \protect\ref{SECCOMPUT})}
{Appendix: Experimental results for the min--max truncated estimator (Section 3.3)}}


%
\begin{sidewaystable}
\textwidth=\textheight
\tablewidth=\textwidth
\tabcolsep=2pt
\caption{Comparison of the min--max truncated estimator $\hf$ with the
ordinary least squares estimator $\hfols$ for the mixture noise (see
Section~\protect\ref{secnoise}) with $\rho=0.1$ and $p=0.005$. In
parenthesis,
the $95\%$-confidence intervals for the estimated quantities}
\label{tabb01}
{\fontsize{7.5pt}{10.5pt}\selectfont{
\begin{tabular*}{\tablewidth}{@{\extracolsep{\fill}}l@{\hspace*{-2pt}}ccccccc@{}}
\hline
& \textbf{Nb of}  & \textbf{Nb of iter. with}
& \textbf{Nb of iter. with}  &
&  & $\bolds{\E R[(\hfols
)|\hf\!\neq\!\hfols]}$
& $\bolds{\E[R(\hf)|\hf\!\neq\!\hfols]}$\\
& \multicolumn{1}{c}{\textbf{iterations}} & $\bolds{R(\hf)\!\neq\!R(\hfols)}$
& $\bolds{R(\hf)\!<\!R(\hfols)}$
& $\bolds{\E R(\hfols)\!-\!R(f^*)}$ & \multicolumn{1}{c}{$\bolds{\E R(\hf)\!-\!R(f^*)}$} &
\multicolumn{1}{l}{$\quad\bolds{-}\,\bolds{R(f^*)}$}
& \multicolumn{1}{l@{}}{$\quad\bolds{-}\,\bolds{R(f^*)}$}\\
\hline
INC $(n\!=\!200,d\!=\!1$)& $1\mbox{,}000$& $419$& $405$& $ 0.567 $ $(\mbox{$\pm$}0.083)$& $ 0.178
$ $(\mbox{$\pm$}0.025)$& $ 1.191 $ $(\mbox{$\pm$}0.178)$& $ 0.262 $ $(\mbox{$\pm$}0.052)$\\
INC $(n\!=\!200,d\!=\!2)$& $1\mbox{,}000$& $506$& $498$& $ 1.055 $ $(\mbox{$\pm$}0.112)$& $ 0.271
$ $(\mbox{$\pm$}0.030)$& $ 1.884 $ $(\mbox{$\pm$}0.193)$& $ 0.334 $ $(\mbox{$\pm$}0.050)$\\
HCC $(n\!=\!200,d\!=\!2)$& $1\mbox{,}000$& $502$& $494$& $ 1.045 $ $(\mbox{$\pm$}0.103)$& $ 0.267
$ $(\mbox{$\pm$}0.024)$& $ 1.866 $ $(\mbox{$\pm$}0.174)$& $ 0.316 $ $(\mbox{$\pm$}0.032)$\\
TS $(n\!=\!200,d\!=\!2)$& $1\mbox{,}000$& $561$& $554$& $ 1.069 $ $(\mbox{$\pm$}0.089)$& $ 0.310
$ $(\mbox{$\pm$}0.027)$& $ 1.720 $ $(\mbox{$\pm$}0.132)$& $ 0.367 $ $(\mbox{$\pm$}0.036)$\\
INC $(n\!=\!1\mbox{,}000,d\!=\!2)$& $1\mbox{,}000$& $402$& $392$& $ 0.204 $ $(\mbox{$\pm$}0.015)$& $ 0.109
$ $(\mbox{$\pm$}0.008)$& $ 0.316 $ $(\mbox{$\pm$}0.029)$& $ 0.081 $ $(\mbox{$\pm$}0.011)$\\
INC $(n\!=\!1\mbox{,}000,d\!=\!10)$& $1\mbox{,}000$& $950$& $946$& $ 1.030 $ $(\mbox{$\pm$}0.041)$& $ 0.228
$ $(\mbox{$\pm$}0.016)$& $ 1.051 $ $(\mbox{$\pm$}0.042)$& $ 0.207 $ $(\mbox{$\pm$}0.014)$\\
HCC $(n\!=\!1\mbox{,}000,d\!=\!10)$& $1\mbox{,}000$& $942$& $942$& $ 0.980 $ $(\mbox{$\pm$}0.038)$& $ 0.222
$ $(\mbox{$\pm$}0.015)$& $ 1.008 $ $(\mbox{$\pm$}0.039)$& $ 0.203 $ $(\mbox{$\pm$}0.015)$\\
TS $(n\!=\!1\mbox{,}000,d\!=\!10)$& $1\mbox{,}000$& $976$& $973$& $ 1.009 $ $(\mbox{$\pm$}0.037)$& $ 0.228
$ $(\mbox{$\pm$}0.017)$& $ 1.018 $ $(\mbox{$\pm$}0.038)$& $ 0.217 $ $(\mbox{$\pm$}0.016)$\\
INC $(n\!=\!2\mbox{,}000,d\!=\!2)$& $1\mbox{,}000$& $209$& $207$& $ 0.104 $ $(\mbox{$\pm$}0.007)$& $ 0.078
$ $(\mbox{$\pm$}0.005)$& $ 0.206 $ $(\mbox{$\pm$}0.021)$& $ 0.082 $ $(\mbox{$\pm$}0.012)$\\
HCC $(n\!=\!2\mbox{,}000,d\!=\!2)$& $1\mbox{,}000$& $184$& $183$& $ 0.099 $ $(\mbox{$\pm$}0.007)$& $ 0.076
$ $(\mbox{$\pm$}0.005)$& $ 0.196 $ $(\mbox{$\pm$}0.023)$& $ 0.070 $ $(\mbox{$\pm$}0.010)$\\
TS $(n\!=\!2\mbox{,}000,d\!=\!2)$& $1\mbox{,}000$& $172$& $171$& $ 0.101 $ $(\mbox{$\pm$}0.007)$& $ 0.080
$ $(\mbox{$\pm$}0.005)$& $ 0.206 $ $(\mbox{$\pm$}0.020)$& $ 0.083 $ $(\mbox{$\pm$}0.012)$\\
INC $(n\!=\!2\mbox{,}000,d\!=\!10)$& $1\mbox{,}000$& $669$& $669$& $ 0.510 $ $(\mbox{$\pm$}0.018)$& $ 0.206
$ $(\mbox{$\pm$}0.012)$& $ 0.572 $ $(\mbox{$\pm$}0.023)$& $ 0.117 $ $(\mbox{$\pm$}0.009)$\\
HCC $(n\!=\!2\mbox{,}000,d\!=\!10)$& $1\mbox{,}000$& $669$& $669$& $ 0.499 $ $(\mbox{$\pm$}0.018)$& $ 0.207
$ $(\mbox{$\pm$}0.013)$& $ 0.561 $ $(\mbox{$\pm$}0.023)$& $ 0.125 $ $(\mbox{$\pm$}0.011)$\\
TS $(n\!=\!2\mbox{,}000,d\!=\!10)$& $1\mbox{,}000$& $754$& $753$& $ 0.516 $ $(\mbox{$\pm$}0.018)$& $ 0.195
$ $(\mbox{$\pm$}0.013)$& $ 0.558 $ $(\mbox{$\pm$}0.022)$& $ 0.131 $ $(\mbox{$\pm$}0.011)$\\
\hline
\end{tabular*}}}
\end{sidewaystable}
%

%
\begin{sidewaystable}
\textwidth=\textheight
\tablewidth=\textwidth
\tabcolsep=2pt
\caption{Comparison of the min--max truncated estimator $\hf$ with the
ordinary least squares estimator $\hfols$ for the mixture noise (see
Section~\protect\ref{secnoise}) with $\rho=0.4$ and $p=0.005$. In
parenthesis,
the $95\%$-confidence intervals for the estimated quantities}
\label{tabb04}
{\fontsize{7.5pt}{10.5pt}\selectfont{
\begin{tabular*}{\tablewidth}{@{\extracolsep{\fill}}l@{\hspace*{-2pt}}ccccccc@{}}
\hline
& \textbf{Nb of}  & \textbf{Nb of iter. with}
& \textbf{Nb of iter. with}  &
&  & $\bolds{\E R[(\hfols
)|\hf\!\neq\!\hfols]}$
& $\bolds{\E[R(\hf)|\hf\!\neq\!\hfols]}$\\
& \multicolumn{1}{c}{\textbf{iterations}} & $\bolds{R(\hf)\!\neq\!R(\hfols)}$
& $\bolds{R(\hf)\!<\!R(\hfols)}$
& $\bolds{\E R(\hfols)\!-\!R(f^*)}$ & \multicolumn{1}{c}{$\bolds{\E R(\hf)\!-\!R(f^*)}$} &
\multicolumn{1}{l}{$\quad\bolds{-}\,\bolds{R(f^*)}$}
& \multicolumn{1}{l@{}}{$\quad\bolds{-}\,\bolds{R(f^*)}$}\\
\hline
INC $(n\!=\!200,d\!=\!1)$& $1\mbox{,}000$& $234$& $211$& $ 0.551 $ $(\mbox{$\pm$}0.063)$& $ 0.409
$ $(\mbox{$\pm$}0.042)$& $ 1.211 $ $(\mbox{$\pm$}0.210)$& $ 0.606 $ $(\mbox{$\pm$}0.110)$\\
INC $(n\!=\!200,d\!=\!2)$& $1\mbox{,}000$& $195$& $186$& $ 1.046 $ $(\mbox{$\pm$}0.088)$& $ 0.788
$ $(\mbox{$\pm$}0.061)$& $ 2.174 $ $(\mbox{$\pm$}0.293)$& $ 0.848 $ $(\mbox{$\pm$}0.118)$\\
HCC $(n\!=\!200,d\!=\!2)$& $1\mbox{,}000$& $222$& $215$& $ 1.028 $ $(\mbox{$\pm$}0.079)$& $ 0.748
$ $(\mbox{$\pm$}0.051)$& $ 2.157 $ $(\mbox{$\pm$}0.243)$& $ 0.897 $ $(\mbox{$\pm$}0.112)$\\
TS $(n\!=\!200,d\!=\!2)$& $1\mbox{,}000$& $291$& $268$& $ 1.053 $ $(\mbox{$\pm$}0.079)$& $ 0.805
$ $(\mbox{$\pm$}0.058)$& $ 1.701 $ $(\mbox{$\pm$}0.186)$& $ 0.851 $ $(\mbox{$\pm$}0.093)$\\
INC $(n\!=\!1\mbox{,}000,d\!=\!2)$& $1\mbox{,}000$& $127$& $117$& $ 0.201 $ $(\mbox{$\pm$}0.013)$& $ 0.181
$ $(\mbox{$\pm$}0.012)$& $ 0.366 $ $(\mbox{$\pm$}0.053)$& $ 0.207 $ $(\mbox{$\pm$}0.035)$\\
INC $(n\!=\!1\mbox{,}000,d\!=\!10)$& $1\mbox{,}000$& $262$& $249$& $ 1.023 $ $(\mbox{$\pm$}0.035)$& $ 0.902
$ $(\mbox{$\pm$}0.030)$& $ 1.238 $ $(\mbox{$\pm$}0.081)$& $ 0.777 $ $(\mbox{$\pm$}0.054)$\\
HCC $(n\!=\!1\mbox{,}000,d\!=\!10)$& $1\mbox{,}000$& $201$& $192$& $ 0.991 $ $(\mbox{$\pm$}0.033)$& $ 0.902
$ $(\mbox{$\pm$}0.031)$& $ 1.235 $ $(\mbox{$\pm$}0.088)$& $ 0.790 $ $(\mbox{$\pm$}0.067)$\\
TS $(n\!=\!1\mbox{,}000,d\!=\!10)$& $1\mbox{,}000$& $171$& $162$& $ 1.009 $ $(\mbox{$\pm$}0.033)$& $ 0.951
$ $(\mbox{$\pm$}0.031)$& $ 1.166 $ $(\mbox{$\pm$}0.098)$& $ 0.825 $ $(\mbox{$\pm$}0.071)$\\
INC $(n\!=\!2\mbox{,}000,d\!=\!2)$& $1\mbox{,}000$& \hphantom{0}$80$& \hphantom{0}$77$& $ 0.105 $ $(\mbox{$\pm$}0.007)$& $ 0.099
$ $(\mbox{$\pm$}0.006)$& $ 0.214 $ $(\mbox{$\pm$}0.042)$& $ 0.135 $ $(\mbox{$\pm$}0.029)$\\
HCC $(n\!=\!2\mbox{,}000,d\!=\!2)$& $1\mbox{,}000$& \hphantom{0}$44$& \hphantom{0}$42$& $ 0.102 $ $(\mbox{$\pm$}0.007)$& $ 0.099
$ $(\mbox{$\pm$}0.007)$& $ 0.187 $ $(\mbox{$\pm$}0.050)$& $ 0.120 $ $(\mbox{$\pm$}0.034)$\\
TS $(n\!=\!2\mbox{,}000,d\!=\!2)$& $1\mbox{,}000$& \hphantom{0}$47$& \hphantom{0}$47$& $ 0.101 $ $(\mbox{$\pm$}0.007)$& $ 0.099
$ $(\mbox{$\pm$}0.007)$& $ 0.147 $ $(\mbox{$\pm$}0.032)$& $ 0.103 $ $(\mbox{$\pm$}0.026)$\\
INC $(n\!=\!2\mbox{,}000,d\!=\!10)$& $1\mbox{,}000$& $116$& $113$& $ 0.511 $ $(\mbox{$\pm$}0.016)$& $ 0.491
$ $(\mbox{$\pm$}0.016)$& $ 0.611 $ $(\mbox{$\pm$}0.052)$& $ 0.437 $ $(\mbox{$\pm$}0.042)$\\
HCC $(n\!=\!2\mbox{,}000,d\!=\!10)$& $1\mbox{,}000$& $110$& $105$& $ 0.500 $ $(\mbox{$\pm$}0.016)$& $ 0.481
$ $(\mbox{$\pm$}0.015)$& $ 0.602 $ $(\mbox{$\pm$}0.056)$& $ 0.430 $ $(\mbox{$\pm$}0.044)$\\
TS $(n\!=\!2\mbox{,}000,d\!=\!10)$& $1\mbox{,}000$& $101$& $98$& $ 0.511 $ $(\mbox{$\pm$}0.016)$& $ 0.499
$ $(\mbox{$\pm$}0.016)$& $ 0.601 $ $(\mbox{$\pm$}0.054)$& $ 0.486 $ $(\mbox{$\pm$}0.051)$\\
\hline
\end{tabular*}}}
\end{sidewaystable}
%

%
\begin{sidewaystable}
\textwidth=\textheight
\tablewidth=\textwidth
\tabcolsep=2pt
\caption{Comparison of the min--max truncated estimator $\hf$ with the
ordinary least squares estimator $\hfols$ with the heavy-tailed noise
(see Section \protect\ref{secnoise})}
\label{taba201}
{\fontsize{7.5pt}{10.5pt}\selectfont{
\begin{tabular*}{\tablewidth}{@{\extracolsep{\fill}}l@{\hspace*{-2pt}}cccd{2.10}d{2.10}d{2.9}c@{}}
\hline
& \textbf{Nb of}  & \textbf{Nb of iter. with}
& \textbf{Nb of iter. with}  &
&  & \multicolumn{1}{c}{$\bolds{\E R[(\hfols
)|\hf\!\neq\!\hfols]}$}
& \multicolumn{1}{c@{}}{$\bolds{\E[R(\hf)|\hf\!\neq\!\hfols]}$}\\
& \multicolumn{1}{c}{\textbf{iterations}} & $\bolds{R(\hf)\!\neq\!R(\hfols)}$
& \multicolumn{1}{c}{$\bolds{R(\hf)\!<\!R(\hfols)}$}
& \multicolumn{1}{c}{$\bolds{\E R(\hfols)\!-\!R(f^*)}$} & \multicolumn{1}{c}{$\bolds{\E R(\hf)\!-\!R(f^*)}$} &
\multicolumn{1}{l}{$\quad\bolds{-}\,\bolds{R(f^*)}$}
& \multicolumn{1}{l@{}}{$\quad\bolds{-}\,\bolds{R(f^*)}$}\\
\hline
INC $(n\!=\!200,d\!=\!1)$& 1\mbox{,}000& 163& 145&  7.72  \mbox{ }(\mbox{$\pm$}3.46)&  3.92  \mbox{ }(\mbox{$\pm$}
0.409)& 30.52  \mbox{ }(\mbox{$\pm$}20.8)&  \hphantom{0}7.20\mbox{ }(\mbox{$\pm$}1.61)\\
INC $(n\!=\!200,d\!=\!2)$& 1\mbox{,}000& 104& \hphantom{0}98& 22.69  \mbox{ }(\mbox{$\pm$}23.14)& 19.18  \mbox{ }(\mbox{$\pm$}
23.09)& 45.36  \mbox{ }(\mbox{$\pm$}14.1)& 11.63\mbox{ }(\mbox{$\pm$}2.19)\\
HCC $(n\!=\!200,d\!=\!2)$& 1\mbox{,}000& 120& 117& 18.16  \mbox{ }(\mbox{$\pm$}12.68)&  8.07  \mbox{ }(\mbox{$\pm$}
0.718)& 99.39  \mbox{ }(\mbox{$\pm$}105)&  15.34\mbox{ }(\mbox{$\pm$}4.41)\\
TS $(n\!=\!200,d\!=\!2)$& 1\mbox{,}000& 110& 105& 43.89  \mbox{ }(\mbox{$\pm$}63.79)& 39.71  \mbox{ }(\mbox{$\pm$}
63.76)& 48.55  \mbox{ }(\mbox{$\pm$}18.4)& 10.59\mbox{ }(\mbox{$\pm$}2.01)\\
INC $(n\!=\!1\mbox{,}000,d\!=\!2)$& 1\mbox{,}000& 104& 100&  3.98  \mbox{ }(\mbox{$\pm$}2.25)&  1.78
 \mbox{ }(\mbox{$\pm$}0.128)& 23.18  \mbox{ }(\mbox{$\pm$}21.3)&  \hphantom{0}2.03\mbox{ }(\mbox{$\pm$}0.56)\\
INC $(n\!=\!1\mbox{,}000,d\!=\!10)$& 1\mbox{,}000& 253& 242&16.36  \mbox{ }(\mbox{$\pm$}5.10)&  7.90
 \mbox{ }(\mbox{$\pm$}0.278)& 41.25  \mbox{ }(\mbox{$\pm$}19.8)&  \hphantom{0}7.81\mbox{ }(\mbox{$\pm$}0.69)\\
HCC $(n\!=\!1\mbox{,}000,d\!=\!10)$& 1\mbox{,}000& 220& 211&13.57  \mbox{ }(\mbox{$\pm$}1.93)&  7.88
 \mbox{ }(\mbox{$\pm$}0.255)& 33.13  \mbox{ }(\mbox{$\pm$}8.2)&  \hphantom{0}7.28\mbox{ }(\mbox{$\pm$}0.59)\\
TS $(n\!=\!1\mbox{,}000,d\!=\!10)$& 1\mbox{,}000& 214& 211& 18.67  \mbox{ }(\mbox{$\pm$}11.62)& 13.79  \mbox{ }(\mbox{$\pm$}
11.52)& 30.34  \mbox{ }(\mbox{$\pm$}7.2)&  \hphantom{0}7.53\mbox{ }(\mbox{$\pm$}0.58)\\
INC $(n\!=\!2\mbox{,}000,d\!=\!2)$& 1\mbox{,}000& 113& 103&  1.56  \mbox{ }(\mbox{$\pm$}0.41)&  0.89
 \mbox{ }(\mbox{$\pm$}0.059)&  6.74  \mbox{ }(\mbox{$\pm$}3.4)&  \hphantom{0}0.86\mbox{ }(\mbox{$\pm$}0.18)\\
HCC $(n\!=\!2\mbox{,}000,d\!=\!2)$& 1\mbox{,}000& 105& \hphantom{0}97&  1.66  \mbox{ }(\mbox{$\pm$}0.43)&  0.95  \mbox{ }(\mbox{$\pm$}
0.062)&  7.87  \mbox{ }(\mbox{$\pm$}3.8)&  \hphantom{0}1.13\mbox{ }(\mbox{$\pm$}0.23)\\
TS $(n\!=\!2\mbox{,}000,d\!=\!2)$& 1\mbox{,}000& 101& \hphantom{0}95&  1.59  \mbox{ }(\mbox{$\pm$}0.64)&  0.88  \mbox{ }(\mbox{$\pm$}
0.058)&  8.03  \mbox{ }(\mbox{$\pm$}6.2)&  \hphantom{0}1.04\mbox{ }(\mbox{$\pm$}0.22)\\
INC $(n\!=\!2\mbox{,}000,d\!=\!10)$& 1\mbox{,}000& 259& 255& 8.77  \mbox{ }(\mbox{$\pm$}4.02)&  4.23
 \mbox{ }(\mbox{$\pm$}0.154)& 21.54  \mbox{ }(\mbox{$\pm$}15.4)&  \hphantom{0}4.03\mbox{ }(\mbox{$\pm$}0.39)\\
HCC $(n\!=\!2\mbox{,}000,d\!=\!10)$& 1\mbox{,}000& 250& 242& 6.98  \mbox{ }(\mbox{$\pm$}1.17)&  4.13
 \mbox{ }(\mbox{$\pm$}0.127)& 15.35  \mbox{ }(\mbox{$\pm$}4.5)&  \hphantom{0}3.94\mbox{ }(\mbox{$\pm$}0.25)\\
TS $(n\!=\!2\mbox{,}000,d\!=\!10)$& 1\mbox{,}000& 238& 233&  8.49  \mbox{ }(\mbox{$\pm$}3.61)&  5.95
 \mbox{ }(\mbox{$\pm$}3.486)& 14.82  \mbox{ }(\mbox{$\pm$}3.8)&  \hphantom{0}4.17\mbox{ }(\mbox{$\pm$}0.30)\\
\hline
\end{tabular*}}}
\end{sidewaystable}
%

%
\begin{sidewaystable}
\textwidth=\textheight
\tablewidth=\textwidth
\tabcolsep=2pt
\caption{Comparison of the min--max truncated estimator $\hf$ with the
ordinary least squares estimator $\hfols$ with the asymmetric
heavy-tailed noise\break (see Section \protect\ref{secnoise})}
\label{taba-201}
{\fontsize{7.5pt}{10.5pt}\selectfont{
\begin{tabular*}{\tablewidth}{@{\extracolsep{\fill}}l@{\hspace*{-2pt}}cd{3.0}d{3.0}d{2.10}d{2.10}d{2.9}c@{}}
\hline
& \textbf{Nb of}  & \multicolumn{1}{c}{\textbf{Nb of iter. with}}
& \multicolumn{1}{c}{\textbf{Nb of iter. with}}  &
&  & \multicolumn{1}{c}{$\bolds{\E R[(\hfols
)|\hf\!\neq\!\hfols]}$}
& \multicolumn{1}{c@{}}{$\bolds{\E[R(\hf)|\hf\!\neq\!\hfols]}$}\\
& \multicolumn{1}{c}{\textbf{iterations}} & \multicolumn{1}{c}{$\bolds{R(\hf)\!\neq\!R(\hfols)}$}
& \multicolumn{1}{c}{$\bolds{R(\hf)\!<\!R(\hfols)}$}
& \multicolumn{1}{c}{$\bolds{\E R(\hfols)\!-\!R(f^*)}$} & \multicolumn{1}{c}{$\bolds{\E R(\hf)\!-\!R(f^*)}$} &
\multicolumn{1}{l}{$\quad\bolds{-}\,\bolds{R(f^*)}$}
& \multicolumn{1}{l@{}}{$\quad\bolds{-}\,\bolds{R(f^*)}$}\\
\hline
INC $(n\!=\!200,d\!=\!1)$& 1\mbox{,}000& 87& 77&  5.49 \mbox{ }(\mbox{$\pm$}3.07)&  3.00 \mbox{ }(\mbox{$\pm$}
0.330)& 35.44 \mbox{ }(\mbox{$\pm$}34.7)&  \hphantom{0}6.85\mbox{ }(\mbox{$\pm$}2.48)\\
INC $(n\!=\!200,d\!=\!2)$& 1\mbox{,}000& 70& 66& 19.25 \mbox{ }(\mbox{$\pm$}23.23)& 17.4 \mbox{ }(\mbox{$\pm$}
23.2)& 37.95 \mbox{ }(\mbox{$\pm$}13.1)& 11.05\mbox{ }(\mbox{$\pm$}2.87)\\
HCC $(n\!=\!200,d\!=\!2)$& 1\mbox{,}000& 67& 66&  7.19 \mbox{ }(\mbox{$\pm$}0.88)&  5.81 \mbox{ }(\mbox{$\pm$}
0.397)& 31.52 \mbox{ }(\mbox{$\pm$}10.5)& 10.87\mbox{ }(\mbox{$\pm$}2.64)\\
TS $(n\!=\!200,d\!=\!2)$& 1\mbox{,}000& 76& 68& 39.80 \mbox{ }(\mbox{$\pm$}64.09)& 37.9 \mbox{ }(\mbox{$\pm$}
64.1)& 34.28 \mbox{ }(\mbox{$\pm$}14.8)&  \hphantom{0}9.21\mbox{ }(\mbox{$\pm$}2.05)\\
INC $(n\!=\!1\mbox{,}000,d\!=\!2)$& 1\mbox{,}000& 101& 92&  2.81 \mbox{ }(\mbox{$\pm$}2.21)&  1.31 \mbox{ }(\mbox{$\pm$}
0.106)& 16.76 \mbox{ }(\mbox{$\pm$}21.8)&  \hphantom{0}1.88\mbox{ }(\mbox{$\pm$}0.69)\\
INC $(n\!=\!1\mbox{,}000,d\!=\!10)$& 1\mbox{,}000& 211& 195&10.71 \mbox{ }(\mbox{$\pm$}4.53)&  5.86
\mbox{ }(\mbox{$\pm$}0.222)& 29.00 \mbox{ }(\mbox{$\pm$}21.3)&  \hphantom{0}6.03\mbox{ }(\mbox{$\pm$}0.71)\\
HCC $(n\!=\!1\mbox{,}000,d\!=\!10)$& 1\mbox{,}000& 197& 185& 8.67 \mbox{ }(\mbox{$\pm$}1.16)&  5.81
\mbox{ }(\mbox{$\pm$}0.177)& 20.31 \mbox{ }(\mbox{$\pm$}5.59)&  \hphantom{0}5.79\mbox{ }(\mbox{$\pm$}0.43)\\
TS $(n\!=\!1\mbox{,}000,d\!=\!10)$& 1\mbox{,}000& 258& 233& 13.62 \mbox{ }(\mbox{$\pm$}11.27)& 11.3 \mbox{ }(\mbox{$\pm$}
11.2)& 14.68 \mbox{ }(\mbox{$\pm$}2.45)&  \hphantom{0}5.60\mbox{ }(\mbox{$\pm$}0.36)\\
INC $(n\!=\!2\mbox{,}000,d\!=\!2)$& 1\mbox{,}000& 106& 92&  1.04 \mbox{ }(\mbox{$\pm$}0.37)&  0.64 \mbox{ }(\mbox{$\pm$}
0.042)&  4.54 \mbox{ }(\mbox{$\pm$}3.45)&  \hphantom{0}0.79\mbox{ }(\mbox{$\pm$}0.16)\\
HCC $(n\!=\!2\mbox{,}000,d\!=\!2)$& 1\mbox{,}000& 99& 90&  0.90 \mbox{ }(\mbox{$\pm$}0.11)&  0.66 \mbox{ }(\mbox{$\pm$}
0.042)&  3.23 \mbox{ }(\mbox{$\pm$}0.93)&  \hphantom{0}0.82\mbox{ }(\mbox{$\pm$}0.16)\\
TS $(n\!=\!2\mbox{,}000,d\!=\!2)$& 1\mbox{,}000& 84& 81&  1.11 \mbox{ }(\mbox{$\pm$}0.66)&  0.60 \mbox{ }(\mbox{$\pm$}
0.042)&  6.80 \mbox{ }(\mbox{$\pm$}7.79)&  \hphantom{0}0.69\mbox{ }(\mbox{$\pm$}0.17)\\
INC $(n\!=\!2\mbox{,}000,d\!=\!10)$& 1\mbox{,}000& 238& 222& 6.32 \mbox{ }(\mbox{$\pm$}4.18)&  3.07
\mbox{ }(\mbox{$\pm$}0.147)& 16.84 \mbox{ }(\mbox{$\pm$}17.5)&  \hphantom{0}3.18\mbox{ }(\mbox{$\pm$}0.51)\\
HCC $(n\!=\!2\mbox{,}000,d\!=\!10)$& 1\mbox{,}000& 221& 203& 4.49 \mbox{ }(\mbox{$\pm$}0.98)&  2.98
\mbox{ }(\mbox{$\pm$}0.091)&  9.76 \mbox{ }(\mbox{$\pm$}4.39)&  \hphantom{0}2.93\mbox{ }(\mbox{$\pm$}0.22)\\
TS $(n\!=\!2\mbox{,}000,d\!=\!10)$& 1\mbox{,}000& 412& 350&  5.93 \mbox{ }(\mbox{$\pm$}3.51)&  4.59
\mbox{ }(\mbox{$\pm$}3.44)&  6.07 \mbox{ }(\mbox{$\pm$}1.76)&  \hphantom{0}2.84\mbox{ }(\mbox{$\pm$}0.16)\\
\hline
\end{tabular*}}}
\end{sidewaystable}
%

%
\begin{sidewaystable}
\textwidth=\textheight
\tablewidth=\textwidth
\tabcolsep=2pt
\caption{Comparison of the min--max truncated estimator $\hf$
with the ordinary least squares estimator $\hfols$ for standard
Gaussian noise}
\label{taba0}
{\fontsize{7.5pt}{10.5pt}\selectfont{
\begin{tabular*}{\tablewidth}{@{\extracolsep{\fill}}l@{\hspace*{-2pt}}ccccccc@{}}
\hline
& \textbf{Nb of}  & \textbf{Nb of iter. with}
& \textbf{Nb of iter. with}  &
&  & \multicolumn{1}{c}{$\bolds{\E R[(\hfols
)|\hf\!\neq\!\hfols]}$}
& \multicolumn{1}{c@{}}{$\bolds{\E[R(\hf)|\hf\!\neq\!\hfols]}$}\\
& \multicolumn{1}{c}{\textbf{iterations}} & $\bolds{R(\hf)\!\neq\!R(\hfols)}$
& \multicolumn{1}{c}{$\bolds{R(\hf)\!<\!R(\hfols)}$}
& \multicolumn{1}{c}{$\bolds{\E R(\hfols)\!-\!R(f^*)}$} & \multicolumn{1}{c}{$\bolds{\E R(\hf)\!-\!R(f^*)}$} &
\multicolumn{1}{l}{$\quad\bolds{-}\,\bolds{R(f^*)}$}
& \multicolumn{1}{l@{}}{$\quad\bolds{-}\,\bolds{R(f^*)}$}\\
\hline
INC $(n\!=\!200,d\!=\!1)$& $1\mbox{,}000$& $20$\hphantom{0}& $8$& $ 0.541 $ $(\mbox{$\pm$}0.048)$& $ 0.541 $ $(\mbox{$\pm$}
0.048)$& $ 0.401 $ $(\mbox{$\pm$}0.168)$& $ 0.397 $ $(\mbox{$\pm$}0.167)$\\
INC $(n\!=\!200,d\!=\!2)$& $1\mbox{,}000$& $1$& $0$& $ 1.051 $ $(\mbox{$\pm$}0.067)$& $ 1.051 $ $(\mbox{$\pm$}
0.067)$& $ 2.566 $& $ 2.757 $\\
HCC $(n\!=\!200,d\!=\!2)$& $1\mbox{,}000$& $1$& $0$& $ 1.051 $ $(\mbox{$\pm$}0.067)$& $ 1.051 $ $(\mbox{$\pm$}
0.067)$& $ 2.566 $& $ 2.757 $\\
TS $(n\!=\!200,d\!=\!2)$& $1\mbox{,}000$& $0$& $0$& $ 1.068 $ $(\mbox{$\pm$}0.067)$& $ 1.068 $ $(\mbox{$\pm$}
0.067)$& --& --\\
INC $(n\!=\!1\mbox{,}000,d\!=\!2)$& $1\mbox{,}000$& $0$& $0$& $ 0.203 $ $(\mbox{$\pm$}0.013)$& $ 0.203 $ $(\mbox{$\pm$}
0.013)$& --& --\\
INC $(n\!=\!1\mbox{,}000,d\!=\!10)$& $1\mbox{,}000$& $0$& $0$& $ 1.023 $ $(\mbox{$\pm$}0.029)$& $ 1.023
$ $(\mbox{$\pm$}0.029)$& --& --\\
HCC $(n\!=\!1\mbox{,}000,d\!=\!10)$& $1\mbox{,}000$& $0$& $0$& $ 1.023 $ $(\mbox{$\pm$}0.029)$& $ 1.023
$ $(\mbox{$\pm$}0.029)$& --& --\\
TS $(n\!=\!1\mbox{,}000,d\!=\!10)$& $1\mbox{,}000$& $0$& $0$& $ 0.997 $ $(\mbox{$\pm$}0.028)$& $ 0.997 $ $(\mbox{$\pm$}
0.028)$& --& --\\
INC $(n\!=\!2\mbox{,}000,d\!=\!2)$& $1\mbox{,}000$& $0$& $0$& $ 0.112 $ $(\mbox{$\pm$}0.007)$& $ 0.112 $ $(\mbox{$\pm$}
0.007)$& --& --\\
HCC $(n\!=\!2\mbox{,}000,d\!=\!2)$& $1\mbox{,}000$& $0$& $0$& $ 0.112 $ $(\mbox{$\pm$}0.007)$& $ 0.112 $ $(\mbox{$\pm$}
0.007)$& --& --\\
TS $(n\!=\!2\mbox{,}000,d\!=\!2)$& $1\mbox{,}000$& $0$& $0$& $ 0.098 $ $(\mbox{$\pm$}0.006)$& $ 0.098 $ $(\mbox{$\pm$}
0.006)$& --& --\\
INC $(n\!=\!2\mbox{,}000,d\!=\!10)$& $1\mbox{,}000$& $0$& $0$& $ 0.517 $ $(\mbox{$\pm$}0.015)$& $ 0.517
$ $(\mbox{$\pm$}0.015)$& --& --\\
HCC $(n\!=\!2\mbox{,}000,d\!=\!10)$& $1\mbox{,}000$& $0$& $0$& $ 0.517 $ $(\mbox{$\pm$}0.015)$& $ 0.517
$ $(\mbox{$\pm$}0.015)$& --& --\\
TS $(n\!=\!2\mbox{,}000,d\!=\!10)$& $1\mbox{,}000$& $0$& $0$& $ 0.501 $ $(\mbox{$\pm$}0.015)$& $ 0.501 $ $(\mbox{$\pm$}
0.015)$& --& --\\
\hline
\end{tabular*}}}
\end{sidewaystable}
%

%

%
\end{appendix}


\newpage

\begin{supplement}[id=suppA]
\stitle{Supplement to
``Robust linear least squares regression''\\}
\slink[doi]{10.1214/11-AOS918SUPP} 
\sdatatype{.pdf}
\sfilename{aos918\_supp.pdf}
\sdescription{The supplementary material provides the proofs of
Theorems \ref{thhfrlam}, \ref{thermom} and \ref{th31}.}
\end{supplement}

%

\printaddresses

\end{document}